\documentclass[reqno,12pt]{amsart}

\usepackage{amsthm,amsfonts,amsxtra,amssymb,amscd}

\usepackage[all]{xy}

\theoremstyle{plain}       
\newtheorem{lem}{Lemma}
\newtheorem{teo}[lem]{Theorem}
\newtheorem{prp}[lem]{Proposition}
\newtheorem{cor}[lem]{Corollary}
\newtheorem{con}[lem]{Conjecture}

\theoremstyle{definition}
\newtheorem{dfn}[lem]{Definition}
\newtheorem{oss}[lem]{Remark}

\theoremstyle{remark}

\newcommand{\mC}  {{\mathbb C}}  
\newcommand{\mR}  {{\mathbb R}}

\newcommand{\mi}  {{\imath}}
\newcommand{\mj}  {{\jmath}}
\newcommand{\mZ}  {\mathbb Z}
\newcommand{\mN}  {{\mathbb N}}

\newcommand{\mP}  {\text{${\mathbb P}$}}
\newcommand{\calB}  {\text{$\mathcal {B}$}}
\newcommand{\calC}  {\mathcal C}

\newcommand{\calF}  {\text{$\mathcal {F}$}}
\newcommand{\calG}  {\text{$\mathcal {G}$}}
\newcommand{\calH}  {\text{$\mathcal {H}$}}
\newcommand{\calI}  {\text{$\mathcal {I}$}}

\newcommand{\calK}  {\mathcal K}

\newcommand{\calM}  {\mathcal M}

\newcommand{\calQ}  {\text{$\mathcal {Q}$}}
\newcommand{\calR}  {\text{$\mathcal {R}$}}
\newcommand{\calS}  {\text{$\mathcal {S}$}}

\newcommand{\grg}  {\text{$\gamma$}}    
\newcommand{\grd}  {\text{$\delta$}}
\newcommand{\grG}  {\text{$\Gamma$}}   
\newcommand{\grD}  {\text{$\Delta$}}
\newcommand{\gro}  {\omega}
\newcommand{\grO}  {\text{$\Omega$}}
\newcommand{\gra}  {\text{$\alpha$}}
\newcommand{\grl}  {\text{$\lambda$}}
\newcommand{\grL}  {\text{$\Lambda$}}
\newcommand{\grm}  {\text{$\mu$}}
\newcommand{\grs}  {\text{$\sigma$}}
\newcommand{\grb}  {\text{$\beta$}}
\newcommand{\gre}  {\varepsilon}
\newcommand{\grf}  {\text{$\varphi$}}
\newcommand{\grz}  {\text{$\zeta$}}
\newcommand{\gog}  {\text{$\mathfrak g$}} 
\newcommand{\goh}  {\text{$\mathfrak h$}}

\newcommand{\gou}  {\text{$\mathfrak u$}}

\newcommand{\goz}  {\text{$\mathfrak z$}}
\newcommand{\goL}  {\text{$\mathfrak L$}}
\newcommand{\goM}  {\text{$\mathfrak M$}}
\newcommand{\goF}  {\mathfrak F}
\newcommand{\goS}  {\mathfrak S}
\newcommand{\goZ}  {\mathfrak Z}

\newcommand{\msec}    {m^{\prime \prime}}

\newcommand{\ssec}    {s^{\prime \prime}}

\newcommand{\vsec}    {v^{\prime \prime}}

\newcommand{\grlsec}  {\grl^{\prime \prime }}
\newcommand{\mter}    {m^{\prime \prime \prime}}
\newcommand{\ster}    {s^{\prime \prime \prime}}
\newcommand{\vter}    {v^{\prime \prime \prime}}

\newcommand{\grlter}  {\grl^{\prime \prime \prime}}
\renewcommand{\sec}   {^{\prime \prime}}
\newcommand{\ter}   {^{\prime \prime \prime}}

\newcommand{\pu}        {\wp}

\newcommand{\isocan}    {\simeq}
\newcommand{\vuoto}     {\varnothing}

\newcommand{\cech}      {\spcheck}

\newcommand{\coinc}     {\equiv}
\newcommand{\circa}	{\sim}
\newcommand{\comp}      {\!\circ\!}
\newcommand{\vuotosq}   {\square}

\newcommand{\im}         {\text{\rm i}} 
\newcommand{\Id}         {\mathrm{Id}}

\newcommand{\rank}	 {\operatorname{rank}}

\newcommand{\tc}         {\, : \,}
\newcommand{\st}         {\, : \,}

\newcommand{\Tr}         {\operatorname{Tr}}    
\renewcommand{\Re}         {\operatorname{Re}}    
\renewcommand{\Im}         {\operatorname{Im}}    
\newcommand{\then}       {\Rightarrow}
\newcommand{\Hom}        {\text{Hom}}

\newcommand{\incluso}    {\hookrightarrow}

\newcommand{\mfor}       {\text{ for }}
\newcommand{\mand}       {\text{ and }}
\newcommand{\mif}        {\text{ if }}
\newcommand{\msuchthat}   {\text{ such that }}

\newcommand{\wt}[1]      { {\widetilde {#1} } }
\newcommand{\wbar}[1]    { {\overline  {#1} } }

\newcommand{\lra}        {\longrightarrow}

\newcommand{\pairing} {< \: ,  >}
\newcommand{\bra} {< \! \!}
\newcommand{\ket} {\! \! >}

\begin{document}

\title{A remark on quiver varieties and Weyl groups}
\author{Andrea Maffei}
\begin{abstract}
In this paper we define an action of the Weyl group
on the quiver varieties $M_{m,\grl}(d,v)$ with generic $(m,\grl)$. To do it we describe a set of 
generators of the projective ring of a quiver variety. We also prove connectness for the smooth quiver variety
$M(d,v)$ and normality for $M_0(d,v)$ in the case of a quiver of finite type and $d-v$ a regular weight.
\end{abstract}

\maketitle

In \cite{Na1,Na2} Nakajima defined quiver varieties and show how to use them to give a geometric
construction of integrable representation of Kac-Moody algebras. Luckily these varieties can be 
used also to give a geometric construction of representations of Weyl groups. 
In \cite{Lu:Q4}, Lusztig constructed a representation of the Weyl group 
on the homology of quiver varieties.
His construction is similar to the construction of Springer representations.
In \cite{Na2},
Nakajima gave an construction of isomorphism $\Phi_{\grs,\zeta}(d,v):\goM_{\zeta}(d,v)
\lra \goM_{\grs \zeta}(d,\grs(v-d)+d)$ in the case of a quiver of finite type. 
His construction was analytic and relies on a description of quiver varieties 
as moduli spaces of instantons on ALE spaces. 
The main result of this paper is a direct and algebraic construction of 
these  isomorphism which works for a general 
quiver without simple loops. 
To do it we also describe a set of generators of the algebra of covariant functions.

The paper is organized as follows. In the first section we fix the notation and we give 
the definition of a quiver variety: $M_{m,\grl}(d,v)$ where $m, \grl$ are two parameter,
$d$ is a weight of the algebra associated to the quiver and $v$ an element of the root lattice.
We are interested to quiver varieties as algebraic varieties but
to explain one of the applications we need to give also the hyperK\"ahler construction of a quiver variety.
We use a result of Migliorini \cite{Migliorini} 
to explain the connection between the two constructions.

Algebraic quiver varieties are defined as the $\mathbf{Proj}$ scheme of a ring of covariants.
In the second section we describe a set of generators of this ring. In a special case which is
not directly to Nakajima's quiver varieties we are also able to give a more precise results and 
to describe a basis of the vector space of $\chi$-covariants functions.

In the third section we use this description to generalize a construction of Lusztig
\cite{Lu:Q4}. Namely for any element of  the Weyl group we construct an isomorphism
$\Phi_{\grs}$ between $M_{m,\grl}(d,v)$ and $M_{\grs m,\grs\grl}(d,\grs(v-d)+d)$ if $m,\grl$
are generic.

In the fourth section, following Nakajima \cite{Na1}, we show how to 
use the action constructed in section 3 
(and the connection between the hyperK\"ahler construction 
and the algebraic construction) to describe an action of the Weyl group on 
the homology of a class of quiver varieties. This action is different from the one
constructed by Lusztig in \cite{Lu:Q4}.

In the fifth section we give a result which reduce the study  of 
geometric and algebraic properties of quiver varieties $M_{0,0}(d,v)$ 
to the case $d-v$ dominant. 

In the sixth section we prove the normality of the quiver 
variety $M_0(d,v)$ and the connectdness of $M(d,v)$ in the 
case of a quiver of finite type and $d-v$ a regular weight.

I wish to thank Ilaria Damiani for many usefull discussions.


\section{Notations and definitions}
In this section we give the definition of quiver varieties. Except
some minor change all definition are due to Nakajima \cite{Na1,Na2}.


\subsection{The graph}
Let $(I,H)$ be a finite oriented graph: 
$I$ is the set of vertices that we suppose of cardinality $n$, 
$H$ the set of arrows and the orientation is given by the two maps 
$$  h \longmapsto h_0 \mand h \longmapsto h_1 $$ from $H$ to $I$. 
We suppose also that:
\begin{enumerate}
  \item $ \forall \, h\in H \quad h_0 \neq h_1$,
  \item an involution $ h\mapsto \bar h $ of $H$ without fixed points and 
 satisfying ${\bar h}_0 = h_1$ is fixed,
  \item a map $\gre : H \lra \{-1,1\}$ is given such that $\gre(\bar h) = - \gre(h)$.
	We define $\grO = \{ h \in H \st \gre(h)= 1\}$ and ${\wbar {\grO}} =
	 \{ h \in H \st \gre(h)= -1\}$.
\end{enumerate}
Observe that given a symmetric graph without loops is always possible to 
define $\gre$ and an involution $\bar {\;}$ as above.


 \subsection{The Cartan matrix and the Weyl group}
Let $A$ be the matrix whose entries are the numbers
$$a_{ij}= card\{ h \in H \, : \, h_0=i \mand h_1=j \}.$$
We define a generalized symmetric Cartan matrix by $C=2I-A$.
Following \cite{Lusztig:QG} an $X,Y$-regular root datum 
$(I,X,X \cech ,\pairing)$ with Cartan matrix equal to $C$ 
is defined in the following way:
\begin{enumerate}
	\item $X \cech$ and $X$ are finetely generated free abelian groups,
	\item $\pairing: X \times X \cech  \lra \mZ$ is a perfect
		bilinear pairing,
	\item two linearly independent sets $\Pi= \{\gra_i \tc i \in I \}
		\subset X$ and 
		 $\Pi \cech =\{ \gra_i\cech \tc i \in I \}
		\subset X \cech$ are fixed and we set $Q = \langle \Pi \rangle$
		and $Q \cech = \langle \Pi \cech \rangle$,
	\item $\bra \gra_i \, ,  \gra_j \cech \ket = c_{ij}$,
	\item (nonstandard) $\rank X = \rank X \cech = 2n - \rank C$,
	\item (nonstandard) a linearly independent set 
		$\{ \gro_i \tc i\in I \}$ of $X$ such that 
		$\bra \gro_i , \gra_j \cech \ket = \grd_{ij}$ is fixed.
\end{enumerate}
Once $C$ is given it is easy to construct a data as above.
We call $\goh$ the complexification of $X \cech$ and we observe that through 
the bilinear pairing $\pairing$ we can identify $\goh^*$ 
with the complexification of $X$. We observe also that 
the triple $(\goh,\Pi,\Pi \cech)$ is a realization  
of the Cartan matrix $C$ (\cite{Kac} pg.1).

The Weyl group $W$ attached to $C$ is defined as the subgroup
of $Aut(X) \subset GL(\goh^*)$ generated by the reflections 
\begin{equation}\label{eq:generatoriW}
s_i: x \longmapsto x - \! \bra x, \gra_i \cech \ket \gra_i.
\end{equation}
Observe that the dual action is given by 
$s_i( y)= y - \! \bra  \gra_i, y \ket \gra_i\cech$ and that 
the lattices $Q$ and $Q\cech$ are stable for these actions.
So the annihilator $\overset{\circ}{ Q\cech} =\{x \in X \tc \bra x,y\ket \, = 0
\: \forall y \in Q \cech\}$ is also stable by $W$ and we can consider the 
action of $W$on the lattice $P=X \, / \, \overset{\circ}{ Q\cech} \isocan 
\Hom_{\mZ}(Q,\mZ)$ and we call $x \mapsto 
\wbar x$ the projection from $X$ to $P$. We observe also that this projection 
is an isomorphism from the lattice $\wt P$, that is not $W$-stable, spanned by 
$\{ \gro_i \tc i \in I\}$ and $P$. Finally we observe that 
$$\wbar {\gra}_i = \sum _{j \in I} c_{ij} \wbar {\gro}_j.$$


\subsection{$d, v$ and the space of all matrices}
For the exposition it will be usefull to identify the set $I$ with the set 
of integers $\{ 1, \dots, n \}$.

Let $d=(d_1,\dots,d_n)$ and $v=(v_1,\dots,v_n)$ 
be two $n$-tuples of integers.
We also think of $d, v$ as elements of $X$ in the following way:
\begin{equation}\label{dvX}
   d = \sum_{i\in I} d_i {\gro}_i \; \mand \; 
v= \sum _{i \in I} v_i {\gra}_i\:;
\end{equation}
and through this identification we define an action of $W$ on $v$.
We define also $v\cech = \sum_{i \in I} v_i \gra_i \cech \in Q \cech$.
Once $d,v$ are fixed we fix complex vector spaces $D_i$ and $V_i$ of 
dimensions $d_i$ and $v_i$ and we define the following spaces of maps:
\begin{subequations} \label{eq:decompS}
\begin{align}
	S_{\grO}(d,v) &= \bigoplus _{i \in I} \Hom(D_i,V_i) \oplus \bigoplus 
			_{h \in \grO} \Hom(V_{h_0},V_{h_1}),\\
	S_{\wbar \grO}(d,v) &=\bigoplus _{i \in I} \Hom(V_i,D_i) \oplus \bigoplus 
			_{h \in \wbar \grO} \Hom(V_{h_0},V_{h_1}), \\
	S(d,v)&= S_{\grO}(d,v) \oplus S_{\wbar \grO}(d,v) .
\end{align}
\end{subequations}
More often, when it will not be ambiguous we will write $S_{\grO}, 
S_{\wbar \grO}$ and $S$ instead of $S_{\grO}(d,v), S_{\wbar \grO}(d,v)$
and $S(d,v)$.

For each $h \in H$ (resp. $i\in I$) we define the projection $B_h$ (resp.
$\grg_i$ and $\grd_i$) from $S$ to $\Hom(V_{h_0},V_{h_1})$ (resp. $\Hom(D_i,V_i)$ 
and $\Hom(V_i,D_i)$) with respect to the decomposition described in \eqref{eq:decompS}.

When an element $s$ of $S$ is fixed we will often write $B_h$ (resp. $\grg_i$,
$\grd_i$) instead of $B_h(s)$ (resp. $\grg_i(s)$ and $\grd_i(s)$).
We will also use 
$\grg$ for $(\grg_1 , \dots , \grg_n)$, $\grd$ for 
$(\grd_1 , \dots , \grd_n)$ and $B$ for $(B_h)_{h \in H}$ and often 
we will write an element of $S$ as a triple $(B,\grg,\grd)$.

Once $D_i,V_i$ and an element $s$ of $S$ are fixed we define also:
\begin{subequations}\label{defTab}
\begin{align}
	T_i &= D_i \oplus \bigoplus _{h \tc h_1 =i} V_{h_0} \, , \label{defT}\\
	a_i &= a_i(s) = (\grd_i(s), (B_{\bar h}(s))_{h \tc h_1 =i}): 
		V_i \lra T_i \, ,\label{defa}\\
	b_i &= b_i(s) = (\grg_i(s), (\gre(h)B_{ h}(s))_{h \tc h_1 =i}): 
		T_i \lra V_i \, .\label{defb}
\end{align}
\end{subequations}

We will identify the dual of space of the 
$\mC$-linear maps $\Hom (E,F)$
between two finite dimensional vector spaces with $\Hom(F,E)$ through the pairing
$\bra \grf , \psi \ket = \Tr ( \grf \comp \psi )$. So we can describe $S$ also 
as $S_{\grO} \oplus S_{\grO}^* = T^* S_{\grO}$ and we observe that a natural 
symplectic structure $\gro$ is defined over $S$ by 
$$
\gro((s_{\grO},s_{\wbar \grO}),(t_{\grO},t_{\wbar \grO})) = 
\bra s_{\grO},t_{\wbar \grO}\ket -\bra t_{\grO},s_{\wbar \grO}\ket.
$$ 


\subsection{Hermitian structure on $S$} 
\label{hyperKahlerstructure}

We suppose now that the spaces $D_i$, $V_i$ are endowed with hermitian metrics.
So we can speak of the adjoint  $\grf^*$ of a linear map between these spaces, and
we have a positive definite hermitian structure $h$ on $S$ with explicit formula:
\begin{align} 
	h( (B,\grg,\grd),(\wt B,\wt \grg,\wt \grd)) &= \sum _{h\in H}
		 \Tr (B_h {\wt B}_{ h}^*) + \sum _{i \in I } 
			\Tr (\grg_i {\wt \grg}_i^* + {\wt \grd}_i^* \grd_i) 
			\label{defh}	\\
		&= \sum_{i \in I} \Tr(a_i {\wt a}_i^* + {\wt b}_i^* b_i)  \notag
\end{align}
and an associated real and closed symplectic form 
$\gro_I(s,t)=\Re h (\im s,  t) - \Im h(s,t)$.

\subsection{Group actions and moment maps}\label{groupactions}
We can define an action of the groups $G=GL(V)=\prod  GL(V_i)$ and 
$GL(D)=\prod  GL(D_i)$ on the set $S$ in the following way:
\begin{align}
	g (B_h,\grg_i,\grd_i) &= 
		 (g_{h_1}B_h g_{h_0}^{-1},g_i \grg _i, 
		\grd_i g_i^{-1}) &&\mfor g=(g_i)\in GL(V), \\
	g (B_h,\grg_i,\grd_i) &= (B_h , \grg _i g_i^{-1}, 
		g_i \grd_i ) &&\mfor g=(g_i)\in GL(D).
\end{align}
Observe that these actions commute and that $\gro$ is $GL(V)$ invariant. 
Moreover if $U=U(V)=\prod U(V_i)$ is 
the group of unitary trasformations in $GL(V)$ the real simplectic form
$\gro_I$ is $U(V)$ invariant.

Define
$\grm, \grm_I \colon S \lra \gog=\oplus gl(V_i)$ 
by the following explicit formulas:
\begin{align*}
	\grm_i(B,\grg,\grd)&= \sum_{h \in H \tc h_1=i} \gre(h) B_h B_{\bar h}
 				+ \grg_i \grd_i = b_i a_i\, , \\
	\grm_{I,i}(B,\grg,\grd)&= \frac{\im}{2}
			\left(\sum _{h\in H\tc h_1=i} B_h B_h^* 
		-   B_{\bar h}^* B_{\bar h} + \grg_i \grg_i^* -\grd_i^* \grd_i 
		\right) =\frac{\im}{2} ( b_i b_i^* - a_i^* a_i),
\end{align*}

If we identify $\gog^*=\Hom_{\mC}(\gog,\mC)$ (resp. $\gou^* = \Hom_{\mR} (\gou , \mR)$)
 with $\gog=\oplus gl(V_i)$ (resp.  $\gou$)
through the pairing $\bra (x_i)\, , (y_i) \ket = \sum _i \Tr (x_i y_i)$,
we can observe 
that $\grm$ is a moment map for the action of $G$ on the symplectic 
manifold $(S, \gro)$ and that $\grm_I$ is a moment map 
for the action of $U$ on the symplectic manifold $(S, \gro_I)$. 
It is common to group all these moment maps together and to 
define an hyperK\"ahler moment map 
$$ 
\wt \grm = (\grm_I, \grm): S \lra \gou \oplus \gog = (\mR \oplus \mC)
\otimes_{\mR} \gou.
$$


\subsection{Quiver varieties as  hyperK\"ahler quotients}
\label{definitionquiver}
Let $\grz_i = (\xi_i,\grl_i) \in \mR \oplus \mC$ and 
$\grz = (\grz_1,\dots,\grz_N)$. We define:
$$
\goL_{\grz}(d,v) = \{s \in S \tc \grm (s) -  \grl_i \Id_{V_i} = 0 \,\mand\, 
		     \grm_{I,i} (s) - \im \xi_i  \Id_{V_i} = 0 \}.
$$
We observe that $\goL_{\grz}(d,v)$ is stable for the action of $U(V)$, so, at least as a 
topological Hausdorf space we can define the \emph{quiver variety of type} 
$\grz$ as
$$
\goM_{\grz} (d,v) = \goL_{\grz}(d,v) / U(V).
$$
It will be convenient to define also $\goM_{\grz} (d,v)=\vuoto$
if $d,v \in \mZ^n $ and there exists $i$ such that $v_i <0$ or  $d_i <0$ 
for some $i$. We call $\goZ = \mR^n \oplus \mC^n$ and we observe that we can idenyify it to 
$(\mR \oplus \mC) \otimes _{\mZ}P$ through:
\begin{equation}\label{zetaRCP}
(\xi_1,\dots,\xi_n,\grl_1,\dots,\grl_n)
\longleftrightarrow \sum _{i \in I} (\xi_i , \grl_i) \wbar{\gro}_i.
\end{equation}
In particular we consider an action of the Weyl group $W$ on $\goZ$
through this identification.

\begin{oss}\label{oss:centrogoZ}
There is a surjective map from 
$\goZ$ to $Z_U(\gou) \oplus Z_G(\gog)$:
$$
(\xi_1,\dots,\xi_n,\grl_1,\dots,\grl_n) \lra
\sum _{i \in I} ( \im \xi_i,\grl_i) \Id_{V_i}
$$ 
Observe that $\goL_{\zeta}$ is the fiber of $\wt \mu$ over the image of $\zeta$ in
$Z_U(\gou) \oplus Z_G(\gog)$. 
\end{oss}

\begin{oss}\label{oss:riduzioneanonzero}
If $v , d \geq 0$ define: $I^* = \{ i \in I \st v_i \neq 0\}$, 
$H^* = \{ h \in H \st h_0,h_1 \in I^*\}$, $\gre^* = \gre\bigr|_{H^*}$,
$v^* = (v_i)_{i \in I^*}$, $d^* =(d_i )_{i \in I^*}$, and 
$\zeta^*= (\grz_i)_{i \in I^*}$
then it is clear that 
$$
\goL_ {\grz^*}(d^*,v^*) \isocan \goL_{\grz}(d,v) \;\mand\;
\goM_ {\grz^*}(d^*,v^*) \isocan \goM_{\grz}(d,v).
$$
\end{oss}

Except for the last equivalence which is trivial in our case the following is a general
well known fact (\cite{GIT} ch. 8).
 
\begin{lem} \label{regolaritamu}
Let $s \in \goL_{\zeta}$ then 
\begin{gather*}
d{\wt {\grm}}_s \text{ is surjective} \iff 
d\grm_s \text{ is surjective} \iff 
 d\grm_I \text{ is surjective } \iff \\
\iff \dim \mathrm{Stab}_G\{s\}= 0
\iff \mathrm{Stab}_G\{s\}=\{1_G\}
\end{gather*}
\end{lem}

\begin{dfn}
If $u \in \mZ^n = Q\cech$ $A\subset Q\cech$ we define
$$
\calH_u=\{\grz =(\xi,\grl)\in \goZ \tc \bra \xi \, , u \cech \ket= \bra \grl \, , u \cech \ket=0 \}
\; \mand \;
\calH_A= \bigcup_{u\in A} H_u.
$$
Let now $U_v = \{  u \in \mN^n-\{0\} \msuchthat$ $0 \leq u_i \leq v_i \}$ and
$\calH = \calH_{U_v}$. $\calH$ is a union of a finite number of 
real subspace of $\goZ$ of codimension $3$.  
\end{dfn}

\begin{lem}[Nakajima, \cite{Na1}] \label{lem:regmu}
If $\grz \in \goZ - \calH $ and $\wt {\grm}(s) = \grz$ then 
$\mathrm{Stab}_G\{s\}={1_G}$.
\end{lem}

As a consequence of the above lemma and general results on 
on hyperK\"ahler manifolds (for example \cite{Hit} or \cite{HKLR})
we obtain the following corollary.
\begin{cor}\label{dimensioniquiver}
If $\grz \in \goz - \calH$ then if it is not empty
$\goM_{\grz} (d,v)$ is a smooth 
hyperK\"ahler manifold of real dimension
$2 \bra 2d-v , v \cech \ket$.  
\end{cor}

\subsection{Geometric invariant theory and moment map} \label{ssec:GIT}
In this section we explain the relation between 
the moment map and the GIT quotient proved 
by Kempf, Ness \cite{KeNe}, Kirwan \cite{GIT} and others.
To be more precise we need a generalization of their results in the case 
of an action on an affine variety
proved by Migliorini \cite{Migliorini}.

Let $X$ be an affine variety over $\mC$ and $G$ a reductive group acting on 
$X$. We can assume that $X$ is a closed subvariety of a vector space $V$ 
where $G$ acts linearly. Let $h$ be an hermitian form on $V$ invariant
by the action of a maximal compact group $U$ of $G$ and define a real
$U$-invariant symplectic form on $V$ by
$$\eta (x,y) = \Re h(\im x,y).$$
Then we can define a moment map $\nu : V \lra \gou^*=\Hom_{\mR}(\gou,\mR) $:
$$\bra \nu(x), u \ket = \tfrac{1}{2} \eta ( u \cdot x , x).$$
We observe that the real symplectic form $\eta$
resctricted to a complex submanifold 
is always non degenerate and that $\mu$ restricted to the non singular locus of 
$X$ is a moment map for the action of $U$ on $X$.

Now let $\chi$ be a multiplicative character of $G$. We observe 
that for all $g \in U$ we have $|\chi(g)|=1$ so 
$\im \,d\chi : u \lra \mR$. In particular we can think to 
$\im\,d\chi$ as an element of $\gou ^*$. Morover we observe that it
is invariant by the dual adjoint action, hence it makes sense
to consider the quotient:
$$\goM = \nu^{-1}(i\,d\chi) / U.$$
As we saw our variety are a particular case of this construction.

\medskip

On the other side we can consider the GIT quotient. 
Let us remind the definition. If $\grf$ is 
a character of $G$ we consider the 
line bundle $L_{\grf} = V \times \mC$ on $V$ with the following 
$G$-linearization:
$$
g(x,z)= (g\cdot x,\grf(g)z).
$$
An invariant section of $L_{\grf}$ is determined by an algebraic function
$f:V\lra \mC$ such that $f(gx)=\grf(g)f(x)$ for all
$g \in G$ and $x \in V$. We use the same symbol $L_{\grf}$ 
also for the restriction of $L_{\grf}$ to $X$. 

Given a rational 
action of $G$ on $\mC$-vector space $A$ we define 
\begin{gather*}
A_{\grf,n} = \{ a \in A \tc g\cdot a = \grf^{-n}(g) a \text{ for all } g \in G\},\\
A_{\grf} = \bigoplus _{n=0}^{\infty} A_{\grf,n} \quad 
\text{ as a graded vector space.}
\end{gather*}
Hence we have that $H^0(X,L_{\grf})^G = \mC[X]_{\grf,1}$.
We observe that if $I$ is the ideal of algebraic function on $V$ 
vanishing on $X$ then 
$$H^0(X,L_{\grf})^G = \frac{H^0(V,L_{\grf})^G}{I_{\grf,1}}.$$
This last fact can be proved easily for example averaging a 
$\grf$ equivariant function $f$ on $X$ in the following way:
$$\tilde f (v) = \int_U \grf^{-1}(u) f(u\cdot v) \, du .$$
\begin{dfn}
A point $x$ of $X$ is said to be $\chi$-\emph{semistable} if
there exist $n>0$ and $f \in H^0(X,L_{\chi}^{\otimes n})^G$ such that
$f(x) \neq 0$. We observe that by the remark above
a point of $X$ is $\chi$-semistable if and only if 
is $\chi$-semistable as a point of $V$. We call $X^{ss}_{\chi}$ 
(resp. $V^{ss}_{\chi}$) the open subset 
of $\chi$-semistable points of $X$ (resp. $V$).
\end{dfn}
 
\begin{prp}[ \cite{GIT,Newstead}]
There exists a good quotient of $X^{ss}_{\chi}$ by the action of $G$ and 
we have that
$$ X^{ss}_{\chi} //G = \mathrm{Proj}\, \mC[X]_{\chi}.$$
Moreover $\mathrm{Proj} \,\mC[X]_{\chi}$ is a finetely generated
$\mC$-algebra and a natural projective map 
$$\pi : X^{ss}_{\chi} //G \lra X//G = \mathrm{Spec} \, \mC[X]^G$$
is defined.
\end{prp}

In the case of $\chi \coinc 1$ the following fact is well known:
$$\mathrm{Proj}\, \mC[X]_{\chi}= \mathrm{Spec}\, \mC[X]^G = 
\nu^{-1}(0)/U.$$
The following result is less well known, and its proof requires some 
adjustment of the classical proof for the case $\chi \coinc 1$ (see
for example an appendix of
\cite{Migliorini} or \cite{TatoPhD} ).

\begin{prp}[Migliorini, \cite{Migliorini}]\label{prp:KNM1}
Let $x \in X$ then 
$$
\exists g \in G \st \nu(gx) = \im d\chi \iff
Gx \text{ is a closed orbit in } X^{ss}_{\chi}.
$$
\end{prp}
\begin{prp}[Migliorini, \cite{Migliorini}]\label{prp:KNM2}
The inclusion $\nu^{-1}(\im d\chi) \subset X^{ss}_{\chi}$ induces 
an homeomorphism
$$
\nu^{-1}(\im d\chi) / U \isocan  X^{ss}_{\chi}//G . 
$$
\end{prp}


\subsection{Quiver varieties as algebraic varieties}
If $m= (m _1,\dots, m _n) \in \mZ^N$ we define a character 
$\chi_m$ of $G_v$ by 
$\chi_m = \prod_{i \in I} \det^{m_i}_{GL(V_I)}$.

If $\grl = (\grl_1,\dots,\grl_n ) \in \mC^n$ and $m=(m_1,\dots , m_n) \in \mZ^n$
then we define the varieties:
\begin{align*}
\grL_{\grl}(d,v) &= \{ s \in S \st \grm_i(s) - \grl_i \Id_{V_i} =0 \text{ for all }i\}, \\
\grL_{m,\grl}(d,v) &= \{ s \in \grL_{\grl}(d,v) \tc  s \text{ is } \chi_m-\text{semistable}\}.
\end{align*}
and the associeted \emph{quiver varieties}
\begin{align*}
M_{m,\grl}(d,v) &= \grL_{m,\grl}(d,v) /\!/ G_v \, \mand \, \\
M^0_{\grl}(d,v) &= M_{0,\grl}(d,v) = \grL_{\grl}(d,v) /\!/G_v 
\end{align*}
We call  $p_{m, \grl}^{d,v} \colon \grL_{m,\grl}(d,v) \lra M_{m,\grl}(d,v)$
the quotient map.
Observe that the inclusion $\grL_{m,\grl}(d,v) \subset  \grL_{\grl}(d,v)$
induces a projective morphism 
$$
 \pi_{m,\grl}^{d,v} : M_{m,\grl}(d,v) \lra M^0_{\grl}(d,v).
$$

Finally it will be convenient to define $M_{m,\grl}(d,v)=\vuoto$ if 
$d,v \in \mZ^n$ and $v_i <0$ or $d_i<0$ for some $i$.

We identify $\mZ^n$ with $P$ and $Z=\mC^n$ with $\mC \otimes_{\mZ} P$ through
$$
(m_1,\dots,m_n) \longrightarrow \sum m_i \wbar \omega_i
\;\mand\;
(\grl_1,\dots,\grl_n) \longrightarrow \sum \grl_i \wbar \omega_i.
$$

\begin{oss}
As in \ref{oss:centrogoZ} we have a surjective map from  $Z$ to  
$Z(\gog)$ and $\grL_{\grl}(d,v)$ is the fiber of $\mu$ over the image of $\grl$ in 
$Z(\gog)$.
\end{oss}
\begin{oss} Remark \ref{oss:riduzioneanonzero} holds without changes also in this case.
\end{oss}
\begin{oss}Observe that $P \oplus Z \subset \goZ$. Observe also 
that the map $m \lra \chi_m$ define a surjective morphism from $P$ to $\Hom(G_v,\mC^*)$
and that the following diagram commute:
$$
\xymatrix{
P \ar@{^{(}->}[d] \ar[rr]& & \Hom(G,\mC^*) \ar@{}[r]|>>>>>{\ni} \ar[d]& \chi \ar[d]\\
\mR^n \ar[rr] && {Z(\gou) \isocan (\gou^*)^U} \ar@{}[r]|>>>>{\ni} & {\im d\chi} 
}
$$
In particular we can apply \ref{prp:KNM2} to the action of $G_v$ 
on $\grL_{\grl}(d,v)$ and we obtain: 
$$
  \goM_{(m,\grl)}(d,v) \isocan M_{m,\grl}(d,v).
$$
\end{oss}

\begin{prp}\label{prp:1.14}
Let $(m,\grl) \notin \calH$ and $s \in \grL_{m,\grl}(d,v)$
then $\mathrm{Stab}_{G_v}(s) = \{1\}$.
\end{prp}
\begin{proof}
As we have already claimed it is enough to prove
$\dim  \mathrm{Stab}_G(s)=\{1\}$. We know that there is a good quotient of
$\grL_{m,\grl}(d,v)$ so it is enough to prove that any closed orbit has maximal dimension.
By Proposition \ref{prp:KNM1} if $G_vs$ is closed in $\grL_{m,\grl}(d,v)$
then there exists $g \in G_v$ such that $\mu_I(gs) = \im d \chi$. 
The thesis follows now form $(m,\grl) \notin \calH$ and 
Lemma \ref{lem:regmu}.
\end{proof}  

\begin{cor}
If $(m,\grl) \notin \calH$ and $M_{m,\grl}(d,v) \neq \vuoto$ then it is a smooth 
algebraic variety of dimension $\bra v\cech , 2d-v \ket$.
\end{cor}


\subsection{Path algebra and $b$-path algebra}
To describe functions on quiver varieties we need some notation about 
the path algebra.

\begin{dfn} \label{defpath}

A \emph{path} $\gra$ in our graph is a sequence
$h^{(m)}\dots h^{(1)}$ such that $h^{(i)} \in H$ and $h^{(i)}_1 = h^{(i+1)}_0$
for $i =1,\dots,m-1$.  We define also $\gra_0=h^{(1)}_0$, $\gra_1= h^{(m)}_1$
and we say that  the length of $\gra$ is $m$. If $\gra_0=\gra_1$ we say that
$\gra$ is a closed path. 
We consider also the empty paths $\vuoto_i$ for $i\in I$ and we define 
$(\vuoto_i)_0=(\vuoto_i)_1=i$.
The product of path is defined in the obvious way.

A $b$-\emph{path} 
$[\grb]$ in our graph is a sequence  $[ i_{m+1}^{r_{m+1}}
\gra^{(m)}i_{m}^{r_m}\dots\gra^{(1)}i_1^{r_1} ]$, that we write between square 
brackets
such that $i_j \in I$, $\gra^{(j)}$ are $B$-path, $r_j \in \mN$ and 
$\gra^{(j)}_0 = i_j $ and $\gra^{(j)}_1= i_{j+1} $ for $j=1,\dots,m$. We
consider also the ``empty'' 
$b$-paths indiced by elements of $I$: $[\vuoto_i]$.
We define $[\grb]_0 = i_1$, $[\grb]_1 = i_{m+1}$ and 
$[\vuoto_i]_0=[\vuoto_i]_1=i$.
The length of $[\grb]$ is 
$\sum_{j=1}^{m+1} r_j +\sum_{j=1}^{m} length(\gra^{j})$ 
and the product of $b$-paths is defined in the obvious way:
$$
[\grb]\cdot[\grb^{\prime}]=
\begin{cases}
0     &\text{if } [\grb^{\prime}]_1 \neq [\grb]_0 \\
{[\grb \grb^{\prime}]}  &\text{if } [\grb^{\prime}]_1 = [\grb]_0 =i
\end{cases}
$$

Given a path $\gra= h^{(m)}\dots h^{(1)}$ and a $b$-path 
$\grb =[ i_{m+1}^{r_{m+1}} \gra^{(m)}\dots$  $\dots i_1^{r_1} ]$
we define an evaluation of $\gra$ and $\grb$ on $S$
in the following way: if $s=(B,\grg,\grd) \in S $ then
\begin{align*}
 \vuoto_i(s)& = 0 \in \Hom(V_i,V_i)\;\mand\;  [\vuoto_i](s) = 0 \in \Hom(V_i,V_i),\\
\gra(s) &= B_{h^{(m)}} \comp \dots \comp B_{h^{(1)}} \in 
		\Hom (V_{\gra_0},V_{\gra_1}),\\
\grb(s) &= (\grg_{i_{m+1}}\comp\grd_{i_{m+1}})^{r_{m+1}} 	
	\comp \gra^{(m)}(s) \comp (\grg_{i_{m}}\comp
\grd_{i_{m}})^{r_{m}}  	\comp \dots \comp \\
	&\quad \circ \dots \circ \gra^{(1)}(s)   \comp 
	(\grg_{i_{1}}\comp \grd_{i_{1}})^{r_{1}} \in 
	\Hom (D_{\grb_0},D_{\grb_1}).
\end{align*}

The path algebra 
$\calR$ is  the vector 
space spanned by paths with the product induced by the product 
of paths. 
If $i,j \in I$ 
we say that an element in $\calR$ is of type $(i,j)$ if it is in the linear span
of the paths with source in $i$ and target in $j$.

The $b$-path algebra 
$\calQ$ is  the vector 
space spanned by $b$-paths with the product induced by the product 
of $b$-paths described above. 
If $i,j \in I$ 
we say that an element in $\calR$ is of type $(i,j)$ if it is in the linear span
of the $b$-paths with source in $i$ and target in $j$.
\end{dfn}
\begin{oss}
We observe that the evaluation on $S$ is a morphism of algebra from $\calR$
to the algebra defined by the morphisms of the category of vector spaces.
Moreover if $f$ is of type $(i,j)$ we observe that
$f(s) \in \Hom (V_i,V_j)$.
\end{oss}


\section{Generators of the projective ring of a quiver variety}
In this section we want to describe a set of generators of the graded ring
$\mC[S] _{\chi}$ and by consequence of the projective ring of a quiver
variety $\mC[\grL_{\grl}]_{\chi}$. More precisely we will give a set of
generators as $\mC[S]^G$-module of its $l$-homogeneous part: $\mC[S] _{\chi,l}$.
This result is a generalization of the one obtained by Lusztig in the case of
invariants: $\chi \coinc 1$. First of all recall his result.

\begin{teo}[Lusztig, \cite{Lu:Q3} theorem 1.3]\label{generatoriinvarianti}
The ring $\mC[S]^G$ is generated by the polynomials:
$$
 s \longmapsto \Tr\left( \gra(s) \right) \;\mand\;
 s \longmapsto \grf\left(\grd_{\grb_1}(s)\grb(s)\grg_{\grb_0}(s)\right) 
$$
for $\gra$ a closed path, $\grb$ a path and   
$\grf \in \left(\Hom(D_{\grb_0},D_{\grb_1}) \right)^*$.
\end{teo}

\subsubsection{Determinants} \label{sssec:det}
To describe our result we do first some general remark.
Forget for a moment our quiver, and suppose to have a finite set of finite
dimensional vector spaces $X_1,\dots,X_k$ of dimensions $u_1,\dots,u_k$ and a
pair of nonnegative integers $(m_i^+, m_i^-)$ for each of them. Finally let
$m^+, m^-$ two nonnegative integers such that
$$ N=\sum_{i=1}^k m_i^+ u_i + m^+ = \sum_{i=1}^k m_i^- u_i + m^-,$$
and two vector spaces $M^+$ and $M^-$ of dimension $m^+,m^-$.
Construct the vector spaces:
$$
Y = \bigoplus _{i=1}^k \mC^{m_i^-} \otimes X_i \oplus M^-, \qquad
Z = \bigoplus _{i=1}^k \mC^{m_i^+} \otimes X_i \oplus M^+
$$
and observe that $\dim Y = \dim Z = N$. Define an action of the general linear
group $GL(X_i)$ of $X_i$ on $Y$ by 
$$
g_i\cdot(\sum_{j=1}^k v_j \otimes x_j + m) = \sum_{j\neq i} v_j \otimes x_j
 + m+ v_i \otimes g_i x_i,
$$
and also a similar action on $Z$. Hence the vector space $\Hom(Y,Z)$ acquires a
natural structure of $G_X = \prod_{i=1}^k GL(X_i)$ module.
If we choose an isomorphism $\grs$ between
$\Hom(\bigwedge^NY,\bigwedge^NZ)$ and $\mC$ we can define a function $det$ on 
$\Hom(Y,Z)$ by
$$
det ( A) = \grs \left( \bigwedge^nA \right).
$$
For simplicity we do not emphasize the role of $\grs$ on this definition, so
strictly speaking, $\det$ is a function defined only up to a nontrivial
constant factor. We observe also that 

$
\bigwedge^nY \isocan \left({\bigwedge^{u_i}X_1}\right)^{\otimes m_1^-}
\otimes \dots \otimes \left({\bigwedge^{u_k}X_k}\right)^{\otimes m_k^-} \otimes
\bigwedge^{m^-}M^- 
$

\noindent (and similarly for $Z$) so an isomorphism $\grs$ is determined if 
we choose orientations, or basis, of $X_j, M^+,M^-$. Finally observe that for
any $g = (g_j) \in G_X$ we have 
$$\det( g \cdot A) = \prod_{i=1}^k (det_{GL(X_i)}(g_i))^{m_i^+ - m_i^-} \;
\det(A).$$

\subsubsection{Description of generators}\label{sssec:generators}
We go back now to our quiver and we describe a set of covariant polynomials on
$S$. Any  character $\chi$ of the group $G_v=GL(V)$ is of the
form $\chi=\chi_m = \prod_{i \in I} det^{m_i}_{GL(V_i)}$. We fix such a character
and we define
\begin{align*}
I^+ &= \{ i \in I \tc  m _i >0 \} \mand  \wt m_i = m _i  \mif i\in I^+,\\ 
I^0 &= \{ i \in I \tc  m _i =0 \} \mand  \wt m_i = 0         \mif i\in I^0,\\
I^- &= \{ i \in I \tc  m _i <0 \} \mand  \wt m_i = - m _i \mif i\in I^-.
\end{align*}
We use now the construction explained in \ref{sssec:det} in the case $X_i= V_i$
and $m_i^+ - m_i^- = m_i$. We choose ordered sets 
$A=(a_1, \dots,a_{m^-}) \subset (\bigcup D_i)^{m^-}$ and 
$B=(b_1, \dots, b_{m^+}) \subset (\bigcup D_i^*)^{m^+}$ 
and we define a function $I:A,B \lra I$ by 
$a \in D_{I(a)}$,  $b \in D^*_{I(b)}$. In the framework described above it is
then possible to set  $M^- = \bigoplus_{i=1}^{m^-} \mC_{a_i}$ and
$M^+ = \bigoplus_{i=1}^{m^+} \mC_{b_i}$. In particular we have 
$$
  Y = \bigoplus_{i\in I} \bigoplus_{h=1}^{m_i^-} V_i^{(h)}  \oplus
\bigoplus_{i=1}^{m^-} \mC_{a_i}, \qquad
  Z = \bigoplus_{i\in I}
\bigoplus_{k=1}^{m_i^+} V_i^{[k]}  \oplus \bigoplus_{i=1}^{m^+} \mC_{b_i}
$$
where $V_i^{(l)}, V_i^{[l]}$ are isomorphic copies of $V_i$.
We choose now elements of the $b$-path algebra 
as follows: 
\begin{enumerate}
\item for any $i,j \in I$ and for any  $1 \leq h \leq m^-_i$,  
$1 \leq k \leq m^+_j$ we choose an element $\gra^{i,h}_{j,k}$ 
of the $b$-path algebra of type $(i,j)$, 
\item for any $i \in I$, $1 \leq h \leq m^-_i$ and for any 
$1\leq l \leq m^+$ we choose an element  $\gra^{i,h}_l$ 
of the $b$-path algebra of type
$(i,I(b_l))$, 
\item for any $1 \leq l \leq m^-$ and for any $j \in I$, $1 \leq k \leq m^+_j$ 
we choose an element $\gra^l_{j,k}$ of the $b$-path algebra
of type $(I(a_l),j)$, 
\item for any $1 \leq l \leq m^-$ and for any $1 \leq l'\leq m^-$ 
we  choose an element $\gra^l_{l'}$ of the $b$-path algebra  of type
$(I(a_l),I(b_{l'}))$. 
\end{enumerate}
We call such a data 
$\Delta = (\{ (m_i^+,m_i^-) \}_{i\in I} ,(m^+,m^-),A,B,
\gra^{i,h}_{j,k},\gra^{i,h}_l,\gra^l_{j,k},\gra^l_{l'})$ 
a $\chi$-data and 
we attach to it a 
$\chi$-covariant function on $S$:
$$
f_{\Delta} (s) = \det \left( \Psi_{\Delta} (s)\right)
$$
where $\Psi_{\Delta} $ is a linear map from $Y$ to $Z$ defined by 
\begin{align*}
[\Psi_{\Delta}] ^{V_i^{(h)}}_{V_j^{[k]}} (s) &=  \gra^{i,h}_{j,k}(s), \\
[\Psi_{\Delta}] ^{V_i^{(h)}}_{\mC_{b_l}} (s) &=   
	b_l \comp \grd_{I(b_l)}\comp \gra^{i,h}_l (s) , \\ 
[\Psi_{\Delta}] ^{\mC_{a_l}}_{V_j^{[k]}} (s) &= 
	\gra^l_{j,k} (s) \comp \grg_{I(a_l)} \bigr|_{\mC a_l}, \\
[\Psi_{\Delta}] ^{\mC_{a_l}}_{\mC_{b_{l'}}}(s) &= 
  b_{l'} \comp \grd_{I(b_{l'})} \comp
\gra^l_{l'}(s)\comp\grg_{I(a_l)}\bigr|_{\mC a_l}. 
\end{align*}

The function $f_{\Delta}$ are a set of generators as $\mC[S]^G$-module of 
$\mC[S]_{\chi,1}$, but we will need to define a smaller set of generators.
To define this set we  give a notion of good $\Delta$.

\begin{dfn}
A data $\Delta$ as above is said to be $\chi$-\emph{good} if it satisfies the
following conditions:
\begin{enumerate}
\item $m_i^+ + m_i^- =  \wt m_i$ for all $i \in I$, 
\item $\gra^{l}_{l'} = 0$ for all $l,l'$,
\item $\gra^*_*$ is an element of the path algebra 
(and not just an element of the $b$-path algebra which is obviously bigger),
\item $card\{(j,k) \tc \gra^{i,h}_{j,k} \neq 0 \} + card\{l\tc \gra^{i,h}_l
\neq 0\} \leq v_i$ for all $i,h$,
\item $card\{(i,h) \tc \gra^{i,h}_{j,k} \neq 0 \} + card\{l\tc \gra^l_{j,k}
\neq 0\} \leq v_j$ for all $j,k$,
\item for all $l$ there exists at most one pair $(i,h)$ such that
$\gra^{i,h}_l \neq 0$, 
\item for all $l$ there exists at most one pair $(j,k)$ such that 
$\gra^l_{j,k} \neq 0$.
\end{enumerate}
\end{dfn}

For the applications the only important point will be the first one.

\begin{prp} \label{prp:covarianti}
The set of polynomials $f_{\Delta}$ with $\Delta$ $\chi$-good generates 
$\mC[S]_{\chi,1}$ as a $\mC[S]^{G_v}$-module. 
\end{prp}

\begin{oss}
Prof. Weyman said me that in the case $D=0$ a similar proposition has been proved by him
and  for arbitrary characteristic.
\end{oss}

\subsection{Some remark on the invariant theory of $GL(n)$} \label{ssec:invGL}
If $V$ is a finite dimensional representation of a linearly reductive Lie
group $G$ and $S$ is a simple representation of $S$ we write $V[S]$ for the
$S$-isotipic component of type $S$ of $V$.

We now fix $n$ and we make some remark on the representations of $GL(n)$.
To any  partition of height less or
equal to $n$ we associate an irreducible representation of $GL(n)$ in the usual
way. If we multiply these representations by a power of the inverse of
determinant representation we obtain a complete list of irreducible
representations of $GL(n)$. If $\grl$ is a partition we call $\grl$ the
transpose partition as usual and we define $\grl^{op}=( \grl_1 - \grl_n ,
\grl_1 - \grl_{n-1}, \dots , \grl_1 - \grl_1)$. We call $\grd$ the determinant
representation of $GL(n)$ and $\gre= 1^n$ the associated partition. Finally we
call $V$ the natural representation. 

\begin{lem}\label{lem:invGL1}
\begin{enumerate}
\item ${L_{\grl}}^* = \grd^{-\grl_1} \otimes L_{\grl^{op}}$,
\item $\Hom_{GL(n)}(\grd^m, L_{\grl} \otimes L_{\mu}) = 
\begin{cases}
\mC &\mif \grl = \mu^{op} +(m-\mu_1)\gre, \\
0   &\text{ otherwise,}\\
\end{cases} $
\item $\Hom_{GL(n)}(\grd^m, L_{\grl} \otimes L^*_{\mu}) = 
\begin{cases}
\mC &\mif \grl = \mu +m\gre, \\
0   &\text{ otherwise.}
\end{cases}$
\end{enumerate}
\end{lem}

\begin{proof}
We prove only 2). 
\begin{align*}
 Hom_{GL(n)}(\grd^m, L_{\grl} \otimes L_{\mu}) &=
    Hom_{GL(n)}(\grd^m\otimes {L_{\mu}}^*, L_{\grl} ) \\
 &= Hom_{GL(n)}(\grd^{m-\mu_1} \otimes {L_{\mu^{op}}}, L_{\grl} )
\end{align*}
If $m \geq \mu_1$ the last group is isomorphic to
$Hom_{GL(n)}( L_{\mu^{op}+(m-\mu_1) \gre}, L_{\grl} )$ and 
if $m< \mu_1$ is isomorphic to $Hom_{GL(n)}( {L_{\mu^{op}}}; L_{\grl+ (\mu_1
-m)\gre})$. In any case the thesis follows. 
\end{proof}

We want now to describe $\Hom_{GL(n)}
\left(\grd^m,V^{\otimes i} \otimes 
\left( V^* \right)^{\otimes j}  \right)$.
To do it we will use Schur-duality.
Remind that the irreducible representations of the groups $S_m$ are
parametrized by partitions of $m$ and call $S_{\grl}$ the irreducible
representation associated with $\grl$. Consider now the action of $S_m$ on
$V^{\otimes m}$ given by permuting the factors. This action commute with  the
$GL(n)$ action. Schur duality asserts that the action of the group $S_m \times
GL(n)$ on $V^{\otimes m}$ decomposes in the following way:
$$V^{\otimes m} = 
\bigoplus_{\substack{\grl \vdash m \\ ht(\grl)\leq n}}
S_{\grl} \otimes L_{\grl}.
$$
We describe a set of elements of $\Hom_{GL(n)}
\left(\grd^m,V^{\otimes i} \otimes \left( V^* \right)^{\otimes j}  \right)$.
Let $m$ be a nonnegative integers a choose a permutation $\grs$ of
$\{1,\dots,i+mn\}$.
To $\grs$ we associate  maps: 
\begin{align*}
 \Phi_{\grs} &: \left(V^{\otimes i} \otimes (V^*)^{\otimes i}\right)^{GL(n)}
 \lra
 V^{\otimes i+mn} \otimes (V^*)^{\otimes i} [\grd^m]  \\
 \Psi_{\grs} &: \left(V^{\otimes i} \otimes (V^*)^{\otimes i}\right)^{GL(n)}
 \lra
 V^{\otimes i} \otimes (V^*)^{\otimes i+mn} [\grd^{-m}]
\end{align*}
by
\begin{align*}
\Phi_{\grs} (t \otimes s) &= \grs(o\otimes \dots \otimes o\otimes t) \otimes s
\\ 
\Psi_{\grs} (t \otimes s) &= t \otimes \grs(o^* \otimes \dots \otimes
o^*\otimes s) \end{align*}
where $o $ is a nonzero vector in $\bigwedge^n V$ and $o^*$ is a non zero
vector in  $\bigwedge^n V^*$ and $t \in V^{\otimes i}$, $s \in (V^*)^{\otimes
i} $. 
\begin{lem}\label{lem:invGL2}
1) If $i \neq j+mn$ then 
$$\Hom_{GL(N)}
\left(\grd^m,V^{\otimes i} \otimes 
\left( V^* \right)^{\otimes j}  \right) = 0.$$

2) If $m>0$ then 
$$ V^{\otimes i+mn} \otimes ( V^* )^{\otimes i} [\grd^m] = 
\sum _{\grs} \Im \Phi_{\grs}.$$

3) If $m>0$ then 
 $$ V^{\otimes i} \otimes ( V^* )^{\otimes i+mn} [\grd^{-m}] = 
\sum _{\grs} \Im \Psi_{\grs}.$$
\end{lem}

\begin{proof}
1) follows directly from lemma \ref{lem:invGL1}. 

2) Let $M=V^{\otimes i+mn} \otimes ( V^* )^{\otimes i}[\grd^m]$ and 
       $N = (V^{\otimes i} \otimes ( V^* )^{\otimes i})^G$. $N$ is a
       $S_i \times S_i$ module, $M$  is a $S_{i+mn} \times S_i$- module and  the
maps $\Phi_{\calI,\grs}$ are equivariant with respect the $S_i$ action on  
$(V^* )^{\otimes i}$. In particular it is enough to prove that  if $\grl $ is a
partition of $i$, $M_{\grl}$ is the $S_{\grl}$-isotipic component of $M$
w.r.t. the $S_i$ action and $N_{\grl}$ the  $S_{\grl}$-isotipic component of
$N$ w.r.t. the $S_i$ action on $(V^* )^{\otimes i}$ then
$$
M_{\grl} = \sum _{\grs} \Phi_{\grs}(N_{\grl}).
$$
By point 3 of  lemma\ref{lem:invGL1} we have that 
\begin{align*}  
 M &= \bigoplus_{\substack{\grl \vdash i \\ ht(\grl)\leq n} }
	\left( S_{\grl+m\gre} \otimes L_{\grl+m\gre} \right) \otimes
	\left( S_{\grl} \otimes L_{\grl}             \right)^* [\grd^m] \\
   &= \bigoplus_{\substack{\grl \vdash i \\ ht(\grl)\leq n} }  S_{\grl+m\gre} 
	\otimes S_{\grl} \otimes  \left( \grd^m \otimes \left( L_{\grl} \otimes
	{L_{\grl}}^*\right)^G\right) 
\end{align*}
In particular $M_{\grl} = S_{\grl+m\gre}  \otimes S_{\grl} $ is an irreducible
representation of $S_{i+mn}\times S_i$. Observe   
$\sum _{\grs} \Phi(N_{\grl})$ is a $S_{i+mn}\times S_i$-submodule of $M_{\grl}$ and that it
is clearly nonzero. So $M_{\grl}=\sum _{\grs} \Phi_{\grs}(N_{\grl})$ as
claimed.

The proof of 3) is equal to the previous one. 
\end{proof}
 
We want now to give a slightly different formulation of the lemma above. 
Let $M = V^{\otimes i} \otimes (V^*)^{\otimes j}$ we want to describe 
$ M^*_{\grd^m} = \{ \grf \in M^* \st g \cdot \grf = \grd^{-m}(g) \grf \}$.
Of course this problem is completely equivalent to the previous one. What we
want to do is to reformulate in a more convenient way for our purposes
the description of a set of generators of $M^*_{\grd^m}$. 
Let $m \geq 0$ and choose $\calI = \{ I_1,\dots,I_m \} $ a
collection of $m$ disjoint subsets of $\{1, \dots ,i+mn\}$ of cardinality
$n$. Let  $I_j = \{ i_{j1}< \dots < i_{jn} \}$ and 
$\{1, \dots ,i+mn\} - \bigcup \calI = \{ j_1 <\dots <j_i \}$.  
To $\calI$ and to a permutation $\grs \in S_i$ we associate elements
$$
 \phi_{\calI,\grs} \in \left(V^{\otimes i+mn} \otimes (V^*)^{\otimes
i}\right)^*_{\grd^m} 
\; \mand \;
\psi_{\calI,\grs} \in \left(V^{\otimes i} \otimes (V^*)^{\otimes
i+mn}\right)^*_{\grd^{-m}}
$$
defined by
\begin{align*}
\phi_{\calI,\grs} (v_1\otimes \dots v_{i+mn} \otimes \grf_1 \dots \grf_i)&=
   \prod_{j=1}^{m} \bra o^*, v_{j1}\wedge \dots \wedge v_{jn}\ket \cdot
   \prod_{h=1}^i \bra v_{j_h}\, , \grf_{\grs_h} \ket  \\ 
\psi_{\calI,\grs}(v_1\otimes \dots v_{i} \otimes \grf_1 \dots
\grf_{i+mn})&= \prod_{j=1}^{m} \bra o, \grf_{j1}\wedge \dots \wedge
\grf_{jn}\ket \cdot \prod_{h=1}^i \bra v_{\grs_h}\, , \grf_{j_h} \ket 
\end{align*} where $o $ is a nonzero vector in $\bigwedge^n V$ and $o^*$ is a
non zero vector in  $\bigwedge^n V^*$.

\begin{lem}\label{lem:invGL3}
1) If $i \neq j+mn$ then 
$ \left(V^{\otimes i} \otimes ( V^* )^{\otimes j} \right)^*_{\grd^{m}} = 0.$

2) If $m\geq 0$ then 
$ \left(V^{\otimes i+mn} \otimes ( V^* )^{\otimes i} \right)^*_{\grd^{m}} \,$ 
 is generated by the functions $ \phi_{\calI,\grs}.$

3) If $m\geq 0$ then 
$\left(V^{\otimes i} \otimes ( V^* )^{\otimes i+mn} \right)^*_{\grd^{-m}}\,$ 
is generated by the functions $\psi_{\calI,\grs}.$
\end{lem}
 
\begin{proof}
The proof is clear by the previous lemma.
\end{proof}
 
\subsection{A special case}\label{ssec:specialcase}
In this section we proove a special case of Proposition
 \ref{prp:covarianti} in which we are able to give a more precise result. 
To simplify the exposition of the proof of Proposition
\ref{prp:covarianti} we will also prove another lemma.

Here and in the following we will use polarization.
If $V$ is finite dimensional vector space then we can define a map
$$\pu : (V^{\otimes n})^* \lra S^n(V^*)\subset \mC[V] \; \text{ through }\;
\pu(\grf)(v) = \grf(v\otimes \dots \otimes v).$$

\begin{lem}
$\pu$ is surjective, moreover if $V$ is a finite dimensional representation of
a reductive group $\grG$,  and $\chi$ is a character of $\grG$ then
$$\pu((V^{\otimes n})^*_{\chi}) = S^n(V^*)_{\chi}$$
where $E_{\chi}$ is the isotipic component of type $\chi^{-1}$ of a $G$
module $E$. 
\end{lem}

\begin{lem}
For $i =1,\dots,n$ let $\grG_i$ be a reductive group, $\chi_i$ be a character
of $\grG_i$ and $E_i$ be a f.d.representation of $\grG_i$. Let $\grG = \prod
\grG_i$, then $E= \otimes _i E_i $ is a representation of $\grG$ and $\chi =
\prod \chi_i$ is a character of $\grG$. Then
$$
E^*_{\chi} = (E_1)^*_{\chi_1} \otimes \dots \otimes (E_n)^*_{\chi_n}.
$$
\end{lem}

Let $J^+$, $J^-$ be two sets of indeces and define 
$ \wt J^+= \{0\} \coprod J^+$, $\wt J^- = \{0\} \coprod J^- $ and $J = J^+ \times J^-$, 
$\wt J = \wt J^+\times \wt J^- - \{(0,0)\}$.
For each 
$i \in \wt J^+$  (resp. $j \in \wt J^+$) choose a
vector space $Y_i$ (resp. $X_j$) 
and define
$X = \bigoplus _{j\in  J^-} X_j$ and $Y = \bigoplus _{i\in  J^+} Y_i$.
Consider the group
$$
G_{XY} = \prod _{i \in J^+} GL(Y_i) \times \prod _{j \in J^-} GL(X_j) 
$$
and its character $c = \prod _{i \in J^+} \det_{GL(Y_i)} \times \left(\prod
_{j \in J^-} \det_{GL(X_j)} \right)^{-1}$.

We fix a  matrix $r=(r_{ij}) _{i \in \wt J^+ , j \in \wt J^-}$ 
of integers such that $r_{i0}=1=r_{0j}$ for all $i,j$  and $r_{00} =-1$ 
and we consider the vector spaces: 
$$
 H^{XY} = H   = \bigoplus_{(i,j) \in \wt J} \Hom(X_j, Y_i)^{\oplus
r_{ij}} \;\mand\;  
 H^{XY}_0 = H_0 = \bigoplus_{(i,j) \in J} 
		\!\!\Hom(X_j,Y_i)   
$$
where we adopt the convention $E^{\oplus n}=\mC^n = 0$ if $n<0$. 
When the spaces $X,Y$ will be clear from the context we will write 
$H$ and $H_0$ insted of $H^{XY}$ and $H^{XY}_0$.
We fix a basis $e^{ij}_m$ of $\mC^{r_{ij}}$ so we have a canonical
identification
\begin{equation} \label{eq:Hbasitensor}
H=\bigoplus_{(i,j) \in \wt J} \Hom(X_j,Y_i)\otimes  \mC^{ r_{ij}} .
\end{equation}
We want to study $c$-equivariant polynomials on $H$. 
If we choose two finite dimensional vector spaces $\tilde A, \tilde B $,
linear maps $\gra : \tilde A \lra X_0$ , $\grb: Y_0 \lra \tilde B$, and
elements $\grf_{ij} \in (\mC^{r_{ij}})^*$ for all
$i,j$ then we can define a map 
$\Phi_{\grf,\gra,\grb} : H \lra H_0 \oplus 
\Hom(\tilde A,Y) \oplus \Hom(X,\tilde B) \subset \Hom(X\oplus \tilde A, Y \oplus \tilde B)$ 
by  
\begin{equation}\label{eq:sssec:casoXY:Phi}
\Phi_{\grf,\gra,\grb} (\sum_{ (i,j) \in \wt J} A_{ij} \otimes v_{ij}) = 
\sum _{ (i,j) \in J} \grf_{ij}(v_{ij}) A^{ij} + 
\sum _{ i\in J^+} A^{i0} \comp \gra +
\sum _{ j\in J^-} \grb \comp A^{0j} 
\end{equation}
where $ A_{ij} \otimes v_{ij} \in \Hom (X_j,Y_i) \otimes  \mC^{ r_{ij}}$.

The following is a special version of \ref{prp:covarianti}.

\begin{lem}\label{lem:lemmaprpcovarianti}
$\mC[H]_{c}$ is generated as a vector space by the following functions:
$$
s \longmapsto \det \bigl(\Phi_{\grf,\gra,\grb} (s)\bigr)
$$
where $\Phi_{\grf,\gra,\grb} : H \lra H_0$ is as above.
\end{lem}

\subsubsection{The special case}
We will study an even more special case in which we are able to prove a better
result that I find nice. In the above setting suppose that $X_0=Y_0=0$ and that $r_{ij}=1$ for all
$i,j$.

Define the following set of matrices:
\begin{align*}
\calS_n &= \{ S = (s_{ij}) \in \mN^{J^+ \times J^-}\st \sum_{i,j} s_{ij}=n \} \\
\calS^{XY}   &= \calS = \{ S = (s_{ij}) \in \mN^{J^+ \times J^-}\st \sum_{j} s_{ij}= \dim Y_i 
		\;\forall i\in J^+\;\\
	&\qquad \qquad\;  \mand \; \sum_i s_{ij} = \dim X_j \;\forall j\in J^- \;\}
\end{align*}
As for $H$ we will write $\calS$ when the spaces $X_j,Y_i$ will be clear from the context.
Observe that $\calS = \vuoto$ if $\sum_j \dim X_j \neq \sum_i \dim Y_i$ and that 
if $N=\sum \dim X_j=\sum \dim Y_i$ then $\calS \subset\calS_N$. For each 
$card(J^+) \times card(J^-)$ matirx $S= (s_{ij})$
we consider $\grf_{ij} \in \mC^*$ given by $\grf_ {ij}(\grl) = s_{ij}\grl$ and we define
$$
\Phi_S = \Phi_{\grf,0,0} \;\mand \; f_S =f_S^{XY} =\det (\Phi_S ).
$$
\begin{prp}\label{prp:27}
$\{ f_S \}_{S \in \calS^{XY}} $ is a basis of $\mC[H]_{c}$.
\end{prp}
\begin{proof}
We have to compute $S^n(H^*)_{c}= (S^n(H))^*_{c}$ for all $n$.
For all $S \in \calS_n$ define 
$$
E_S = \bigotimes_{(i,j)\in J^+\times J^-} S^{s_{ij}}\left(\Hom (X_j,Y_i) \right).
$$
Observe that $S^n(H) = \bigoplus _{S\in \calS_n} E_S$ as 
a $G$-module. So $S^n(H)^*_{c} = \bigoplus _{S\in \calS_n} (E_S)^*_{c}$.
Observe now that $E_S $ is a quotient of 
\begin{equation}\label{eq:EStilde}
\tilde E_S = \bigotimes _{(i,j) \in J^+ \times J^-} (X_j^*)^{\otimes s_{ij}} \otimes 
   Y_i^{\otimes s_{ij}}.
\end{equation}
By the lemmas in the previous section we have that
$$
(\tilde E_S)^*_{c} = 
\begin{cases}
0 &\mif S \notin \calS^{XY}, \\
\mC & \mif S \in \calS^{XY}.
\end{cases}
$$
So in particular $(E_S)^*_{c} = 0$ if $S \notin \calS^{XY}$. Hence
$\dim S^n(H)^*_{c} \leq card(\calS^{XY})$.

The function $f_S$ are clearly $c$-equivariant so the only thing that we have to prove is 
that they are linearly independent. To prove it we will prove a generalization of it.

If $i \in J^+$ and $j \in J^-$  let $E_{ij}$ be the $card(J^+) \times card(J^-)$ matrix
with a $1$ in the $(i,j)$ position and $0$ elsewhere.

For each $i\in J^+, j\in J^-, m \in \mN$ and $N\in \mN$ we consider the following sentence
$P_{i,j,m,N}$:
\begin{quote}
If $\sum_j X_j = N = \sum_i Y_i$ then $\{ f_{S + mE_{ij}} \}_{S \in \calS^{XY}}$ is lineraly
independent.
\end{quote}
In the case $m=0$ we call this proposition $P_{0,N}$ since it does not depend on $i,j$ and observe
that $\forall N \,P_{0,N}$ is equivalent to our thesis.

For each $N \in \mN$ we consider also the following sentence $Q_N$:
\begin{quote}
If $\sum_j X_j = N = \sum_i Y_i$ then $P_{i,j,m,N}$ is true for all $i\in J^+$, $j\in J^-$ and
$m\in \mN$.
\end{quote}
\noindent \emph{First remark:} $N=1$ is true.

\noindent \emph{Second remark:} let $\calS_0 ^{XY} = \{ S \in S^{XY} \st s_{ij}=0 $ for all $i\in J^+$ and 
$j\in J^-$ such that $\dim Y_i ,\dim X_j \geq 2$. Observe that $\{f_S\}_{S \in \tilde{\calS} ^{XY}}$ 
is linearly independent.

Now we prove $Q_N$ by induction on $N$.

\noindent \emph{Firts step:} $ Q_{N-1} \then P_{0,N}$. Suppose that there exists $c_S \in \mC$
such that 
$$
\sum_{S \in \calS^{XY}} c_S f_S = 0.
$$
If $\dim X_{j_0},\dim Y_{i_0} \geq 2 $ choose a nonzero element $x_{j_0} \in X_{j_0}$ 
(resp.  $y_{i_0} \in Y_{i_0}$) and an hyperplane $ X' _{j_0} \subset X_{j_0}$
(resp.  $ Y'_{i_0} \subset Y_{i_0}$) such that 
$X_{j_0} = \mC x_{j_0} \oplus  X'_{j_0}$ 
(resp. $Y_{j_0} = \mC y_{i_0} \oplus  Y'_{i_0}$) and define:
\begin{equation}\label{eq:prp27:defXYprimo}
\tilde X_j = 
\begin{cases}
X_j &\mif j \neq j_0 \\
X'_{j_0} &\mif j=j_0
\end{cases}
\;\mand \; 
\tilde Y_i = 
\begin{cases}
Y_i &\mif i \neq i_0 \\
Y'_{i_0} &\mif i=i_0
\end{cases}
\end{equation}
and define $\Psi : H^{\tilde X \tilde Y} \lra H^{XY}$ by
\begin{equation}\label{eq:prp27:defPsi}
\Psi(T)\biggr|_{\tilde X_j} = T \quad\mand\quad
\Psi(T)(x_{j_0}) = y_{j_0}.
\end{equation}
Then we see that 
\begin{gather*}
 0=\sum_{S \in \calS^{XY}} c_S f_S^{XY} (\Psi(T)) = 
 \sum_{S \in \calS^{XY}\st s_{i_0j_0}\neq 0} s_{i_0j_0} c_S f_S^{\tilde X \tilde Y} (T) \\
 = \sum_{S \in \calS^{\tilde X \tilde Y}} (s_{i_0j_0}+1) c_{S+E_{i_0j_0}} 
           f_{S+E_{i_0j_0}}^{\tilde X \tilde Y} (T)
\end{gather*}
By induction $P_{i,j,1,N-1}$ is true for all $i,j$ so we see that $c_S=0$ for all 
$S \in \calS$ such that there exists $i_0,j_0$ such that $s_{i_0,j_0}\geq 1$ and 
$\dim X_{j_0} ,\dim Y_{i_0} \geq 2$.
Now we conclude by the second remark.

\noindent \emph{Second step:} $Q_{N-1} \then P_{i_0,j_0,m,N}$ if
 $\dim X_{j_0},\dim Y_{i_0} \geq 2$ and $m\geq 1$. Suppose that 
$ \sum_{S \in \calS^{XY}} c_S f_{S+mE_{i_0j_0}}^{XY} =0$. We can construct 
$\tilde X_j, \tilde Y_i, \Psi$ as in the first step and we see that
\begin{gather*}
 0=\sum_{S \in \calS^{XY}} c_S f_{S+mE_{i_0j_0}}^{XY} (\Psi(T)) = 
 \sum_{S \in \calS^{XY}} (s_{i_0j_0}+m) c_S f_{S+mE_{i_0j_0}}^{\tilde X \tilde Y} (T) \\
  = \sum_{S \in \calS^{\tilde X \tilde Y}} (s_{i_0j_0}+m+1) c_{S+E_{i_0j_0}}
           f_{S+(m+1)E_{i_0j_0}}^{\tilde X \tilde Y} (T)
\end{gather*}
and by $P_{i_0,j_0,m+1,N-1}$ we deduce $c_S = 0$ for all $S$.

\noindent \emph{Third step:} $Q_{N-1} \then Q_N$. By the previous two step we have only 
to prove $P_{i_0,j_0,m,N}$ for $m\geq 1$ and $\dim X_{j_0} = 1 $ or $\dim Y_{i_0} =1$.
We will suppose $\dim X_{j_0} = 1$, the other case is completely similar.
Suppose that $\sum_{S\in \calS^{XY}} c_S f_{S+mE_{i_0j_0}} =0$.
Set 
$$
\tilde {\calS}_i = \{ S \in \calS^{XY} \st s_{ij_0}=1 \}
$$
and observe that since $\dim X_{j_0}=1$ then $\calS^{XY} = \coprod \tilde {\calS}_i$.
Now choose a non zero vector $x_{j_0} \in X_{j_0}$ and for all $i\in J^+$ choose 
a non zero vector $y_{i} \in Y_i$ and an hyperplane $Y_i' $ of $Y_i$ such that 
$Y_i = \mC y_i \oplus Y_i'$.

Now fix $i_1 \neq i_0$ such that $\dim Y_{i_1} \geq 2 $ and consider 
$\tilde J^+ = J^+$ and $\tilde J^- = J^- - \{j_0\}$. For all $i \in \tilde J^+$ and 
for all $j \in \tilde J^-$ define:
$$
\tilde X_j = X_j,
\;\mand\;
\tilde Y_i = 
\begin{cases}
Y_i & \mif i \neq i_1, \\
Y_{i_1}' &\mif i =i_1. 
\end{cases}
$$
For any $S\in \tilde {\calS}_{i_1}$ we define also $t(S) \in \calS^{\tilde X \tilde Y}$ by 
$t(S)_{ij} = s_{ij}$ for all $i \in  \tilde J^+$, $j \in \tilde J^-$. $S\longmapsto t(S)$ 
is a bijection between $\tilde {\calS}_{i_1}$ and $\calS^{\tilde X \tilde Y}$: we call $t^{-1}$ 
the inverse map. 
Finally we define $\Psi : H^{\tilde X \tilde Y} \lra H^{XY}$ as in the previou step and we observe 
that if $S\in \calS$ then $f_{S+mE_{i_0j_0}} \comp \Psi = 0$ if $S\notin \tilde {\calS}_{i_1}$.
Hence
\begin{gather*}
0= \sum_{S\in \calS^{XY}}c_S f_{S+mE_{i_0j_0}}(\Psi(T))= 
   \sum_{S\in \tilde {\calS}_i} c_S f_{t(S)}(T) \\
 = \sum_{S\in \calS^{\tilde X \tilde Y}} c_{t^{-1}(S)} f_{S}(T) 
\end{gather*}
and applying $P_{0,N-1}$ we obtain $c_S=0$ for all $S \in \tilde {\calS}_{i_1}$ if 
$\dim Y_{i_1} \geq 2$ and $i_1 \neq i_0$.

In a similar way we prove $c_S$ if $S\in \tilde {\calS}_{i_1}$ and $\dim Y_{i_1}=1$ and 
$i_1 \neq i_0$. 

Finally we observe that if $S \in \tilde {\calS} _{i_0}$ then 
$f_{S+mE_{i_0j_0}} = (m+1) f_S$, hence $c_S=0$ follows now from $P_{0,N}$ that 
we already know to be true.
\end{proof}

\subsubsection{Proof of Lemma \ref{lem:lemmaprpcovarianti}}
We study first $(H^{\otimes n})_{c}^*$ and then we apply polarization.
As in the previous section we can decompose $H^{\otimes n}$ in summands of the following form:
\begin{equation}\label{eq:lemdefE}
E = \bigotimes_{(i,j)\in \tilde J} (\Hom(X_j,Y_i) \otimes \mC^{r_{ij}})^{\otimes s_{ij}}=
\bigotimes _{ (i,j) \in \wt J} 
(X_{j}^*)^{\otimes s_{ij}} \otimes  Y_{i}^{\otimes s_{ij}} 
\otimes (\mC^{r_{ij}})^{\otimes s_{ij}}
\end{equation}
where $s_{ij}$ are nonnegative integers such that $\sum _{i,j} s_{ij} = n$. 
Observe that the order of the factors is not important for us since we will aply polarization.

We can describe easily $E^*_{c}$ using the
lemma in the previous section. In particular a necessary and sufficient
condition  for the existence of $c$-covariants is 
$\sum _{i\in \wt J^+} s_{ij} = \dim X_j$ for all $j \in J^-$ and 
 $\sum _{j\in \wt J^+} s_{ij} = \dim Y_i$ for all $i \in J^+$. 
Moreover 
$$
E^*_{c} \isocan
\bigotimes _{ (i,j) \in J} 
\left((\mC^{r_{ij}})^*\right)^{\otimes s_{ij}}
\otimes 
\bigotimes _{j \in J^-} (Y_0^*) ^{\otimes s_{0j}}
\otimes
\bigotimes _{i \in J^+} X_0^{\otimes s_{i0}}
$$

To write explicit formulas we choose an order on
the factors of $E$ for example choosing 
a lexicographic order in $i \in \wt J^+$,
$j \in \wt J^-$ and $1\leq q \leq s_{ij}$:
$$
E= 
\underbrace{
X_{1}^* \otimes  Y_{1} \otimes \mC^{r_{11}}}_{q=1}
\otimes 
\dots
\otimes 
\underbrace{
X_{1}^* \otimes  Y_{1} \otimes \mC^{r_{11}}}_{q=s_{11}}
\otimes 
\underbrace{
X_{1}^* \otimes  Y_{2} \otimes \mC^{r_{12}}}_{q=1}
\otimes 
\cdots
$$
Once we have chosen such an order we can write an element of $E$ as linear
combination elements of the form 
$\otimes _{(i,j,q) \in \wt K} x^{i,j,q} \otimes y^{i,j,q} \otimes v^{i,j,q}$
with 
$x^{i,j,q} \in X_j^* $, $ y^{i,j,q} \in Y_i$ $v^{i,j,q} \in \mC^{r_{ij}}$ and
we setted $\wt K= \{ (i,j,q) \in \wt J \times \mN \st 1\leq q \leq s_{ij} \}$.
We define also $K= \{(i,j,q)  \in \wt K \st (i,j) \in J\}$.
Using this convention if 
\begin{equation}\label{eq:defphi18}
\phi = 
\bigotimes _{(i,j,q) \in \wt K} 
\phi ^{i,j,q}  \in 
\bigotimes _{(i,j,q) \in K} 
(\mC^{r_{ij}})^*
\otimes 
\bigotimes _{(0,j,q) \in \wt K} Y_0^* 
\otimes
\bigotimes _{(i,0,q) \in \wt K} X_0
\end{equation}
the corresponding $c$ equivariant linear function on $E$ is defined on an
element 
$s=\otimes _{(i,j,q) \in \wt K} x^{i,j,q} \otimes y^{i,j,q} \otimes
v^{i,j,q}$   by 
\begin{align*}
\phi(s) &= 
\prod _{i \in J^+} \bra \!\!
\bigwedge_{\substack{\lra \\ (i,j,q)\in \wt K} }\!\!\! y^{i,j,q} , o_i^*\ket  
   \prod _{j \in J^-} \bra \!\!
\bigwedge_{\substack{\lra \\ (i,j,q)\in \wt K} }\!\!\! x^{i,j,q} , o_j \ket  \\
& \quad \prod_{(i,j,q) \in K} \phi^{i,j,q}(v^{i,j,q})
\prod_{(i,0,q) \in \wt K} \phi^{i,0,q}(x^{i,0,q})
\prod_{(0,j,q) \in \wt K} \phi^{0,j,q}(y^{0,j,q})
\end{align*}
Now consider the group  
$$\goS = \goS_1\times \goS_2 \times \goS_3
= \prod _{(i,j)\in J} S_{s_{ij}} 
\times
\prod _{ j \in J^-} S_{s_{0j}}
\times 
\prod _{ i \in J^+} S_{s_{i0}}.$$
This group acts naturally on 
$ \bigotimes _{ (i,j) \in J} 
\left((\mC^{r_{i_1j_1}})^*\right)^{\otimes s_{ij}}
\otimes 
\bigotimes _{j \in J^-} (Y_0^*) ^{\otimes s_{0j}}
\otimes 
\bigotimes _{i \in J^+} X_0^{\otimes s_{i0}} = E^*_{c}$ 
by permuting the factors and we observe that 
$$
\pu \bigl( (\grs_1, \grs_2, \grs_3) \phi\bigr) = 
\gre(\grs_2) \gre(\grs_3) \pu (\phi)
$$
for all  $\phi \in E^*_{c}$ and for all $(\grs_1,\grs_2,\grs_3) \in \goS$. 
So we have that 
$$
\pu (E^*_{c}) = 
\pu \biggl( \bigotimes _{ (i,j) \in J} 
S^{s_{ij}} \bigl((\mC^{r_{ij}})^* \bigr) 
\otimes
\bigotimes_{j\in J-}  \bigwedge ^{s_{0j}} Y_0^* 
\otimes
\bigotimes_{i\in J^+} \bigwedge ^{s_{i0}} X_0
\biggr).
$$
In particular since $S^m(V)$ is spanned by vectors of the form 
$v \otimes \dots \otimes v$, $\pu(E^*_{c})$ is spanned by the functions
$\pu(\phi)$ with $\phi $ of the following special form:
\begin{equation}\label{eq:phispeciale}
\phi = 
\bigotimes _{(i,j) \in J} 
(\phi ^{i,j})^{\otimes s_{ij}} 
\otimes 
\bigotimes_{j \in J^-}  \phi^{0,j,1}\wedge \dots \wedge \phi ^{0,j,s_{0j}}
\otimes
\bigotimes_{i \in J^+}  \phi^{i,0,1}\wedge \dots \wedge \phi ^{i,0,s_{i0}}.
\end{equation}

The lemma now follows from the following claim: 

\noindent {\bf{Claim:}} For each $\phi$ as in \eqref{eq:phispeciale}
$\pu(\phi)$ is a linear combination of the functions $\det(\phi_{\grf,\gra,\grb})$.

We prove the claim as follows: we construct vector spaces $A_i$, $B_j$ and 
$A= \bigoplus_{i\in J^+}A_i$, $B=\bigoplus_{j\in J^-} B_j$ and
$$
\tilde H_0 = H_0 \oplus \bigoplus_{i\in J^+} \Hom(A_i,Y_i) \oplus \bigoplus _{j\in J^-}
\Hom(X_j,B_j) \subset \Hom(X\oplus A, Y\oplus B) = \tilde H.
$$
Observe that on $\tilde H, \tilde H_0$ there is an action of 
$\tilde G = G_{XY}\times G_{AB} = G_{XY} \times \prod_{i\in J^+} Gl(A_i) \times 
\prod_{j\in J^-} Gl(B_j)$ and we call $\tilde c$ the character of $\tilde G$ given by
$ (\prod_j \det_{GL(X_j)} \times \prod_i \det_{GL(A_i)} )^{-1} \times
  \prod_i \det_{GL(Y_i)} \times \prod_j \det_{GL(B_j)}$. 
We have an embedding of $G_{XY}$ in $\tilde G$ such that $\grs^* \tilde c = c$.
Observe also that by Proposition \ref{prp:27} we know that the $\tilde c$-covariants functions
on $H$ are generated by the functionts $\det(\Phi_{\tilde S})$ with 
$S \in \tilde {\calS}$: we put a tilde to 
emphasize  that we have to consider also the components  $\{A_i\}$ and $\{B_j\}$.
Then we construct a $G_{XY}$-equivariant map 
$\rho : H \lra \tilde H_0$ such that
\begin{enumerate}
 \item there exists a $\tilde c$-covariant function $f$ on $\tilde H$ such that 
$\pu(\phi)= f \comp \rho$.
 \item for all $\tilde S \in \tilde {\calS}$ there exists 
$\grf, \gra, \grb$ as in equation \eqref{eq:sssec:casoXY:Phi} such that
$\det (\Phi_{\tilde S}) \comp \rho = \det (\Phi_{\grf,\gra,\grb})$ . 
\end{enumerate}
The claim now follows by Proposition \ref{prp:27}.

For $i\in J^+$ and $j\in J^-$ define
$$
 A_i = \mC^{s_{i0}},  \qquad A = \bigoplus_{i\in J^+} A_i, \qquad
 B_j = \mC^{s_{0j}},  \qquad B = \bigoplus_{j\in J^-} B_j.
$$
Define also $\gra_i : A_i \lra X_0$ and $(\grb_j)^t : B_j^* \lra Y_0^*$ by 
\begin{align*}
 \gra_i(e_l) &= \phi^{i,0,l},   &\mand\quad\gra &= \coprod_{i\in J^+} \gra_i : A \lra X_0 \\
 (\grb_i)^t(e^l) &= \phi^{0,j,l}  &\mand \quad \grb^t &= \prod_{j\in J^-} B_j :B^* \lra Y_0,
\end{align*}
where $e_l$ (resp. $e^l$) is the canonical basis of $\mC^m$ (resp. $(\mC^m)^*$).
We define $\grb_i$ (resp. $\grb$)  as the transpose of $(\grb_i)^t$ (resp. $\grb^t$). Now define 
$\rho^{ij} : \Hom(X_j,Y_i) \otimes \mC^{r_{ij}} \lra \Hom (X_j, Y_i)$,
$\rho^{i0} : \Hom(X_0,Y_i) \lra \Hom (A_i, Y_i)$,
$\rho^{0j} : \Hom(X_j,Y_0) \lra \Hom (X_j, B_j)$
by
$$
\rho^{ij}(T \otimes v) = \phi^{i,j} (v) T, \qquad
\rho^{i0}(T)  = T \comp \gra_i, \qquad
\rho^{0j}(T)  = \grb_j \comp T,
$$
and finally define $\rho = \bigoplus_{i,j\in \tilde J} \rho^{ij}: H \lra \tilde H_0$.
Observe that $\rho$ is $G_{XY}$-equivariant. 

Observe now that $\tilde H _0^{\otimes n}= \bigoplus \tilde E_{\tilde S}$ where 
$\tilde S \in \tilde{\calS}$ and $\tilde E_{\tilde S}$ is defined as in \eqref{eq:EStilde}. 
In particular we choose the following summund of $\tilde H _0^{\otimes n}$:
$$
\tilde E  = \bigotimes _{(i,j) \in J} \Hom(X_j,Y_i) ^{\otimes s_{ij}} \otimes 
\bigotimes_{j \in J^-} \Hom(X_j,B_j)^{\otimes s_{0j}} \otimes
\bigotimes_{i \in J^+} \Hom(A_i,Y_i)^{\otimes s_{i0}}
$$
and we observe that $(\tilde E)^*_{\tilde c} = \mC$. Choose a non zero element 
$\tilde {\phi}  \in (\tilde E)^*_{\tilde c}$ and observe that up to a scalar we have
\begin{equation}\label{eq:phiphitilderho}
\pu_{\tilde H}(\tilde{\phi}) \comp \rho = \pu(\phi).
\end{equation}
To see this choose $\phi $ as in \eqref{eq:phispeciale}, and  bases 
$y^{i}_h$ of $Y_i$, $x^{j}_k$ of $X_j^*$ (and its dual basis $z^{j}_k$ of $X_j$) .
Choose also a bases $\gre^{ij}_m$ of $\mC^{r_{ij}}$ such that 
$\phi^{i,j}(\gre^{ij}_m) = \grd_{m,1}$ and set 
$A^{ij} = \rho^{ij}(s)= \sum_{h,k} a^{ij}_{hk} y^i_h\otimes x^j_k$ for $s \in H$. 
Then 
{\allowdisplaybreaks
\begin{align*}
\pu(\phi)(t) &= 
  \sum _{h \in \calK_Y , k \in \calK_X} \; 
	\prod _{i \in J^+} \bra \!\!
		\bigwedge_{\substack{\lra\\ (i,j,q) \in \wt K}}
		a^{ij}_{h(i,j,q) k(i,j,q)} y^i_{h(i,j,q)} , o^*_i \ket
	\prod _{j \in J^-} \bra \!\!
		\bigwedge_{\substack{\lra\\ (i,j,q) \in \wt K}}
		x^j_{k(i,j,q)} , o_j \ket \\*
	& \qquad
	\prod _{i \in J^+} \bra \!\!
		\bigotimes _{\substack{\lra\\ (i,0,q) \in \wt K}}
		x^0_{k(i,0,q)} , \phi^{i,0,1} \wedge \dots \wedge
		\phi^{i,0,s_{i0}} \ket  \\*
	& \qquad
	\prod _{j \in J^-} \bra \!\!
		\bigotimes _{\substack{\lra\\ (0,j,q) \in \wt K}}
		a^{0j}_{h(0,j,q) k(0,j,q)} y^0_{h(0,j,q)} , 
		\phi^{0,j,1} \wedge \dots \wedge
		\phi^{0,j,s_{i0}} \ket  \\
  	& = \sum _{k \in \calK_X}  
	\prod _{i \in J^+} \bra \!\!
		\bigwedge_{\substack{\lra\\ (i,j,q) \in \wt K}}
		A^{ij} z^j_{k(i,j,q)} , o^*_i \ket
	\prod _{j \in J^-} \bra \!\!
		\bigwedge_{\substack{\lra\\ (i,j,q) \in \wt K}}
		x^j_{k(i,j,q)} , o_j \ket \\*
	& \qquad
	\prod _{i \in J^+} \bra \!\!
		\bigwedge _{\substack{\lra\\ (i,0,q) \in \wt K}}
		x^0_{k(i,0,q)} , \phi^{i,0,1} \wedge \dots \wedge
		\phi^{i,0,s_{i0}} \ket  \\*
	& \qquad
	\prod _{j \in J^-} \bra \!\!
		\bigwedge _{\substack{\lra\\ (0,j,q) \in \wt K}}
		A^{0j} z^j_{k(0,j,q)}, 
		\phi^{0,j,1} \wedge \dots \wedge
		\phi^{0,j,s_{i0}} \ket  
\end{align*}
}
where the indeces are as follows:
\begin{align*}
  \calK_X &= \{ k: \wt K \lra \mN \st 1\leq k(i,j,q) \leq \dim X_j\} \\
  \calK_Y &= \{ h: \wt K \lra \mN \st 1\leq h(i,j,q) \leq \dim Y_i\}.
\end{align*}
The lefthandside in \eqref{eq:phiphitilderho} clearly furnishes the same expression.

Finally if we fix 
$\tilde S = (s_{MN})_{N\in \{X_j\}\cup\{A_i\} \mand M\in \{Y_i\}\cup\{B_j\}} \in \tilde {\calS}$ 
and we choose 
$\grf_{ij}=s_{Y_iX_j} \phi^{i,j}$ and $\gra = \coprod_{i \in J^+} s_{Y_iA_i} \gra_i : A \lra X_0$
and $\grb = \coprod_{j \in J^-} s_{B_jX_j} \grb_j : Y_0 \lra B$ we have 
$$
\det (\Phi_{\tilde S}) \comp \rho = \det (\Phi_{\grf,\gra,\grb}).
$$

\begin{oss}
The basis of $\mC[H]_c$ we have described are different from the polarization of the 
natural basis of $E^*_c$. The relation between the two basis is given by formulas of 
the following types
\begin{enumerate}
\item If 
$A = \begin{pmatrix}a_{11} & a_{12} \\ a_{21} & a_{22} \end{pmatrix}$
and $B = \begin{pmatrix}b_{11} & b_{12} \\ b_{21} & b_{22} \end{pmatrix}$ then
$$
\det \begin{pmatrix}a_{11} & b_{12} \\ a_{21} & b_{22} \end{pmatrix}
+
\det \begin{pmatrix}b_{11} & a_{12} \\ b_{21} & a_{22} \end{pmatrix}
 = \det (A+B) - \det A - \det B.
$$ 
\item If $A = \begin{pmatrix}a_{11} & a_{12} \\ a_{21} & a_{22} \end{pmatrix}$,
$B = \begin{pmatrix}b_{11} & b_{12} \\ b_{21} & b_{22} \end{pmatrix}$,
$C =\begin{pmatrix}c_{11} & c_{12} \\ c_{21} & c_{22} \end{pmatrix}$
and $D =\begin{pmatrix}d_{11} & d_{12} \\ d_{21} & d_{22} \end{pmatrix}$
then 
\begin{gather*}
\det \begin{pmatrix}a_{11} & b_{11} \\ a_{21} & b_{21} \end{pmatrix}
\det \begin{pmatrix}c_{12} & d_{12} \\ c_{22} & d_{22} \end{pmatrix}
-\det \begin{pmatrix}a_{11} & b_{12} \\ a_{21} & b_{22} \end{pmatrix}
\det \begin{pmatrix}c_{12} & d_{11} \\ c_{22} & d_{21} \end{pmatrix} + \\
-\det \begin{pmatrix}a_{12} & b_{11} \\ a_{22} & b_{21} \end{pmatrix}
\det \begin{pmatrix}c_{11} & d_{12} \\ c_{21} & d_{22} \end{pmatrix}
+\det \begin{pmatrix}a_{12} & b_{12} \\ a_{22} & b_{22} \end{pmatrix}
\det \begin{pmatrix}c_{11} & d_{11} \\ c_{21} & d_{21} \end{pmatrix} = \\
= -\det \begin{pmatrix}A & B \\ C & D \end{pmatrix}
+\det \begin{pmatrix} A & 0 \\ 0 & D \end{pmatrix}
+\det \begin{pmatrix} 0 & B \\ C & 0 \end{pmatrix}
\end{gather*}
\end{enumerate}
The first type of formula correspond to the reduction of Lemma 
\ref{lem:lemmaprpcovarianti} to the case $r_{ij}=1$ and $X_0=Y_0=1$.
The second type of formula correspond to the case of Proposition 
\ref{prp:27}.
\end{oss}

\subsection{Proof of Proposition \ref{prp:covarianti}}
Choose basis $\calB_i$ (resp. $\calB_i^*$) of the
vector spaces  $D_i$ and $D_i^*$ and we write our vector space $S(d,v)$ in the
following way: 
$$
S = \bigoplus _{h\in H} V_{h_0}^* \otimes V_{h_1} \oplus 
    \bigoplus _{\substack{i \in I \\ b \in \calB_i}} V_{i,b} \oplus
    \bigoplus _{i \in I, b^* \in \calB_i^*} V_{i,b^*}^* 
$$
where $V_{i,b}$ (resp. $V^*_{i,b^*}$) is an isomorphic  copy of
$V_i$ (resp. $V_i^*$). 
We fix also a character $\chi_m$ and $m_i$, $\wt m _i,m_i^+,m_i^-,I^+,I^-,I^0$ as
in \ref{sssec:generators} and we describe first the $\chi_m$-covariants of 
$S^{\otimes n}$. To do it we observe that we can decompose $S^{\otimes n}$
in the following way:
$$
S^{\otimes n} = \bigoplus _{\ell} E^{(\ell)}_1 \otimes \dots \otimes
E^{(\ell)}_n 
$$
where each $E^{(\ell)}_i$ is a representation of $G$ of one of the following
types: $ V_{h_0}^* \otimes V_{h_1}$, $V_{i,b}$ or $V_{i,b^*}^*$. So it is
enough to compute the $\chi$-covariants of each piece 
$E^{(\ell)}_1 \otimes \dots \otimes E^{(\ell)}_n$. We fix one of them: 
$E = E_1 \otimes \dots \otimes E_n$ and we compute $E^*_{\chi}$.
Let  $I^*$ be a copy of $I$ and  fix an isomorphism $i\longleftrightarrow
i^*$ between the two sets. For each $j=1,\dots,n$ we define  a subset
$\calS_j$ of $I\coprod I^*$ according to the following rule:
$$
\calS_j = 
\begin{cases}
 \{h^*_0, h_1\} &\mif E_j = V_{h_0}^*  \otimes V_{h_1}, \\
 \{i\}		&\mif E_j = V_{i,b}  \text{ for some } b \in \calB_i, \\
 \{i^*\}	 &\mif E_j = V_{i,b^*}^* \text{ for some } b^* \in \calB_i^*.
\end{cases}
$$
Let now be $\calS =\coprod_{j=1}^{n} \calS_j$. An element of $\calS$ can be
thought as a couple $(i,j)$ (or $(i^*,j)$) where $i$ (or $i^*$) is in
$\calS_j$. We consider now a special class of partitions of $\calS$:  a
collection $\goF=\{\calC, \calM^{(l)}_i \mfor i\in I \mand 1\leq l\leq m_i\}$ 
of disjoint subsets  of $2^{\calS}$  is called $m$-\emph{special} if:
\begin{enumerate} 
\item $\bigcup \goF$ is a partition of $\calS$,
\item $\forall \, C \in \calC$ $\, card C =2 $ and $\exists i \in I$,
	 $\calS_{j_1}, \calS_{j_2}$ such that $i\in \calS_{j_1}$, $i^*\in
	\calS_{j_2}$ and $C=\{(i,j_1),(i^*,j_2)\}$,   
\item $\forall M \in \calM^{(l)}_i$ we have $M=\{(i,j)\}$ if $i \in I^+$ 
	and  $M=\{(i^*,j)\}$ if $i \in I^-$,  
\item $card \calM ^{(l)}_i = v_i =\dim V_i$. 
\end{enumerate} 

We can represents a special collection with an enriched graph
whose vertices are the sets $\calS_j$ and 
completed   according with the following rules:
\begin{enumerate}
\item we put an arrow from $\calS_{j_1}$ to $\calS_{j_2}$ if there exists
      $C = \{(i,j_2),(i^*,j_1)\}\in \calC$,
\item we put an indexed circle box $\comp^l_i$ on $\calS_j$ if there exists
      $M = \{(i,j)\} \in \calM^{(l)}_i$
\item we put an indexed square box $\vuotosq^l_i$ on $\calS_j$ if there exists
      $M = \{(i^*,j)\} \in \calM^{(l)}_i$
\item if $E_j$ is of type $V_{i,b}$ or $V_{i,b^*}^*$ then we add the element
	$b$ or $b^*$ at the left of the corresponding vertex.  
\item if $E_j$ is of type $V_{h_0}^* \otimes V_{h_1}$ then we write $h$ at the
	left of the corresponding vertex.
\end{enumerate}
Observe that a vertex can be marked with a circle and a square but that it cannot be marked 
with two circles or two squares.

There is a perfect bijection between $m$-special collection $\goF$ and 
graphs as above such that:
\begin{enumerate}
\item the cardinality of vertexes marked with $\comp^l_i$ is $v_i$ for each $i \in I^+$ and
   $1\leq l \leq m_i$,
\item  the cardinality of vertexes marked with $\vuotosq^l_i$ is $v_i$ for each $i \in I^-$ and
   $1\leq l \leq -m_i$.
\end{enumerate}
We will use the same letter $\goF$ to indicate the collection or the graph.

To a special collection $\goF$ as above we attach a function 
$\phi_{\goF}$ on $E$. We define it by the formula
$$
\phi_{\goF}(e_1 \otimes \dots \otimes e_n) =
\prod_{C \in \calC} \phi_C  \cdot
\prod_{i \in I^+} \prod _{l=1}^{m_i} \bra o_i^* , \bigwedge \calM^{(l)}_i\ket
\cdot
\prod_{i \in I^-} \prod _{l=1}^{m_i} \bra o_i , \bigwedge \calM^{(l)}_i\ket
$$
where $o_i$ is a non zero element in $\bigwedge^{v_i}V_i$, $o_i^*$ 
is a non zero element in $\bigwedge^{v_i}V_i^*$ and
\begin{align*}
e_j &= \begin{cases}
	x^*_j \in V_i^*		& \mif  E_j = V_{i,b}^*,\\
	y_j   \in V_i  		& \mif  E_j = V_{i,b},  \\
	x^*_j \otimes y_j \in V^*_{h_0} \otimes V_{h_1} & \mif
		E_j=V^*_{h_0} \otimes V_{h_1}, 	
	\end{cases} \\
\phi_C &= \bra x^*_{j_1} \,, v_{j_2}\ket \qquad \mif
	C=\{(i^*,j_1),(i,j_2)\}\\  
\bigwedge \calM_i^{(l)} &= y_{j_1}\wedge \dots \wedge y_{j_{v_i}}  \quad \mif
	\calM_i^{(l)} =\{ \{(i,j_1)\} , \dots,\{(i,j_{v_i})\} \} \mand i \in
	I^+ \\ 
\bigwedge \calM_i^{(l)} &= x^*_{j_1}\wedge \dots \wedge x^*_{j_{v_i}}  \quad
	\mif \calM_i^{(l)} =\{ \{(i,j_1)\} , \dots,\{(i,j_{v_i})\} \} \mand i
	\in I^-  
\end{align*}
Finally we extend $\phi_{\goF}$ to all $E$ by linearity.
By the lemma above and the discussion in \ref{ssec:invGL} we deduce easily
the following lemma:
\begin{lem}\label{lem:28}
$E^*_{\chi}$ is generated by the functions $\phi_{\goF}$. 
\end{lem}

Proposition \ref{prp:covarianti} now follows from lemma \ref{lem:28}
and the following claim:
\noindent {\bf claim:} for any special collection
$\goF$  the function $\pu(\phi_{\goF})$ is a $\mC[S]^{G_v}$-linear combination
of the functions $f_{\Delta}$ described in \ref{sssec:generators}.

We consider the connected components of the graph. There are only
five possible types of paths:
\begin{enumerate}
\item closed paths,
\item straight paths leaving from a non boxed vertex and arriving in a non boxed
vertex,
\item straight paths leaving from a non boxed vertex and arriving in a circle boxed
vertex,
\item straight paths leaving from a square boxed vertex and arriving in a non
boxed vertex,
\item straight paths leaving from a square boxed vertex and arriving in a
circle boxed vertex.
\end{enumerate}

Let now $\goF_0$ be the union of the connected components of the first two types
and $\goF_1$ be the union of the remaing components. Observe that
$$
\pu(\phi_{\goF}) = \pu(\phi_{\goF_0}) \pu(\phi_{\goF_1}). 
$$
Observe also that  $\phi_{\goF_0}$ is an invariant function (indeed this part of the graph 
corresponds to the situation studied by Lusztig in \cite{Lu:Q4}).
Since we are interested in generators of $\mC[S]_{\chi,1}$as a $\mC[S]^G$-module, 
we can suppose for simplicity $\goF = \goF_1$.

Observe now that each connected component $\grG$ of the graph of the third type
and with a circle $\comp^l_{i_1}$ at the end, has an initial vertex which is an 
$\calS_j = \{i_0^*\}$ and that is marked with $b \in \calB_{i_0}^*$ on the left.
All the other vertexes of the connected component are of type
$\calS_j = \{ h_0^*,h_1\}$ and they define a path 
$\gra^{\grG} $ such that $\gra^{\grG}_0 = i_0$ and 
$\gra^{\grG}_1 = i_1$. We call $b = b({\grG})$  and $l=L_1(\grG)$.

In the same way we see that:
\begin{enumerate}
\item each connected component $\grG$ of the fourth type determines 
a path $\gra^{\grG}$, $b^*=b^*(\grG) \in \calB^*_{\gra^{\grG}_1}$ and $l=L_0(\grG)$
such that $1 \leq l \leq -m_{\gra^{\grG}_0}$,
\item each connected component $\grG$ of the fifth type determines 
a path $\gra^{\grG}$, $l_0=L_0(\grG)$ and $l_1=L_0(\grG)$ such that 
$1 \leq l_0 \leq -m_{\gra^{\grG}_0}$ and $1 \leq l_1 \leq m_{\gra^{\grG}_1}$.
\end{enumerate}

Now we prove the claim in the following way, we construct $X_j$ and $Y_i$ as in 
\ref{ssec:specialcase}, a groups homomorphism $\grs : G_v \lra G_{XY}$ such that
$\grs^* c =\chi_m$, a $G_v$ equivariant 
map $\rho: S \lra H$, and a $G_{XY}$ $c$-covariant function $f$ on $H$
such that:
\begin{enumerate}
 \item for all $\grf,\gra,\grb$ there exists a $\chi_m$-good data such that
 $\det(\Phi_{\phi,\gra,\grb}) \comp \rho = f_{\Delta}$,
 \item $\pu( \phi_{\goF}) = f \comp \rho$ 
\end{enumerate} 
The claim will clearly  follow.

Set 
\begin{align*}
J^- &= \{ (i,l) \st i \in I^- \mand 1\leq l \leq -m_i \}, \\
J^+ &= \{ (i,l) \st i \in I^+ \mand 1\leq l \leq m_i \}.
\end{align*}
For all $(i,l) \in J^-$  choose $X_{(i,l)} = V_i$ and for each 
$(i,l) \in J^+$  choose $Y_{(i,l)} = V_i$. For each 
$(i_0,l_0) \in J^-$ and for each $(i_1,l_1) \in J^+$ define:
\begin{align*}
r_{(i_0,l_0)(i_1,l_1)} &= card\{\, \text{connected component  $\grG$ of the fifth
type  such } \\
 &\qquad \text{that }\gra^{\grG}_0 = i_0,\; \gra^{\grG}_1 = i_1,\; L_0(\grG) = l_0 \mand
  L_1(\grG) = l_1\}
\end{align*}
We the connected component $\grG$ of the set in the left handside as a basis $e_{\grG}$ of the 
vector space $\mC^{r_{(i_0,l_0)(i_1,l_1)}}$. This basis plays the role 
of the basis $e^{ij}_m$ we used to give the identification in 
\eqref{eq:Hbasitensor}.

For each $\grG$ of the third type choose a one dimensional vector space 
$\mC_{b(\grG)}$ and fix a generator $b_{\grG}$.
For each $\grG$ of the fourth type choose a one dimensional vector space 
$\mC_{b^*(\grG)}$ and fix a generator $b^*_{\grG}$.
\begin{align*}
X_0 &= 
\bigoplus_{\grG \text{ of the third type}} \mC_{b(\grG)} =
\bigoplus_{\grG \text{ of the third type}} \mC b_{\grG}
\\
Y_0 &= 
\bigoplus_{\grG \text{ of the fourth type}} \mC_{b^*(\grG)} =
\bigoplus_{\grG \text{ of the fourth type}} \mC b^*_{\grG}.
\end{align*}

Now for each connected component $\grG$ of the third type define 
$\rho^{\grG}  : S \lra \Hom(\mC_{b(\grG)} , Y_{(\gra^{\grG}_1,L_1(\grG))}) $ by 
$$
s \longmapsto \{ \grl \mapsto \gra^{\grG}(s) \grg_{\gra_0^{\grG}} (b(\grG))\grl\}
$$
In a similar way define $\rho^{\grG}$ if $\grG$ is the fourth or of the fifth
type. Finally define 
$$
\rho : S \lra H \; \text{ by }\;
\rho = \bigoplus_{\grG} \rho^{\grG}.
$$
Define also a group homomorphism $\grs : G_v \lra G_{XY}$ by 
$(\grs(g_i)) _{X_{(i_0,l_0)}} = g_{i_0}$ and 
$(\grs(g_i)) _{Y_{(i_1,l_1)}} = g_{i_1}$, and observe that 
$\rho$ is $G_v$ equivariant.

Now we describe $\phi \in (H^{\otimes \tilde n})^*_c$ 
(in general $\tilde n$ is less or equal to $n$)
such that
\begin{equation}\label{eq:puphigoF}
\pu (\phi_{\goF}) (s) = \pu (\phi)(\rho(s)).
\end{equation}
We describe $\phi$ by giving a summunds $\tilde E$ of $H^{\otimes \tilde n}$
as in \eqref{eq:lemdefE} and $\phi \in \tilde E^*_c$ as in 
\eqref{eq:defphi18}.To define $\tilde E$ we have to define 
$s_{(i_1,l_1)(i_0,l_0)}$, $s_{(i_1,l_1)0}$ and $s_{0(i_0,l_0)}$
for all $(i_1,l_1) \in J^+$ and for all $(i_0,l_0)\in J^-$.
We set
\begin{align*}
s_{(i_1,l_1)(i_0,l_0)} & = r_{(i_1,l_1)(i_0,l_0)} \\
s_{(i_1,l_1)0} & = card \{\, \text{ connected component $\grG$ of the third type} \\
		&\qquad \text{ such that } 
				\gra^{\grG}_1 = i_1 \mand L_1(\grG)=l_1 \} \\ 
s_{0(i_0,l_0)} & = card \{\, \text{ connected component $\grG$ of the fourth type} \\
		&\qquad \text{ such that } 
				\gra^{\grG}_0 = i_0 \mand L_0(\grG)=l_0 \}
\end{align*}
Observe that we can choose a bijection $ q \longleftrightarrow \grG_q$ between 
$\{1,\dots,s_{(i_1,l_1)(i_0,l_0)}\}$ and the set of connected component $\grG$ 
of the fifth type such that 
$\gra^{\grG}_1 = i_1$,$ L_1(\grG)=l_1$,
$\gra^{\grG}_0 = i_0 $ and $ L_0(\grG)=l_0 $. So we can define $\phi^{(i_1,l_1),(i_0,l_0),q}$ by
$$
\phi^{(i_1,l_1),(i_0,l_0),q} (e_{\grG}) = \grd_{\grG,\grG_q}.
$$
Observe also that we can choose a bijection $ q \longleftrightarrow \grG_q$ between 
$\{1,\dots,s_{(i_1,l_1)0}\}$ and the set of connected component $\grG$ 
of the third type such that 
$\gra^{\grG}_1 = i_1$,$ L_1(\grG)=l_1$.
So we can define $\phi^{(i_1,l_1),0,q}$ by
$$
\phi^{(i_1,l_1),0,q} (b_{\grG}) = \grd_{\grG,\grG_q}.
$$
In a similar way define $\phi^{0,(i_0,l_0),q}$.

Up to a sign which depends on our choices and ordering 
equation \eqref{eq:puphigoF} is a tautologically satisfied.

Observe now that by lemma \ref{lem:lemmaprpcovarianti} and linearity 
$\mC[H]_c$ is generated by functions $s\longmapsto \det(\Phi_{\grf,\gra,\grb}(s))$
where $\tilde A,\tilde B,\grf,\gra,\grb$ are as in \eqref{eq:sssec:casoXY:Phi} and 
moreover there exists a basis $e_1,\dots,e_{r_A}$ of $\tilde A$ and a basis 
$\tilde e_1,\dots,\tilde e_{r_B}$ of $\tilde B^*$ 
such that for all $i$ there exist a connected component of the 
third type $\grG^A_i$ such that $\gra(e_i)= b_{\grG^A_i}$ 
and for all $i$ there exists a connected component of the fourth type $\grG^B_i$
such that $\tilde e _i (\grb (b^*_{\grG})) = \grd_{\grG\grG^B_i}$. So it is enough to prove 
that if $\tilde A,\tilde B,\grf,\gra,\grb$ are as above then 
there exists a $\chi_m$-good $\grD$ such that 
$$
\det(\Phi_{\grf,\gra,\grb}) \comp \rho  = f_{\Delta}.
$$
We define 
{\allowdisplaybreaks
\begin{align*}
A & = ( b_{\grG^A_i} )_{i=1,\dots,r_A} \in \left(\bigcup D_i\right)^{\dim \tilde A} \\
B & = ( b^*_{\grG^B_i} )_{i=1,\dots,r_B} \in \left(\bigcup D^*_i\right)^{\dim \tilde B} \\
\gra^{ik}_{jh} & = \sum_{ 1\leq q \leq s_{(j,h)(i,k)}} \phi^{(j,h),(i,k),q}(e_{\grG}) \gra^{\grG}\\
\gra^{ik}_l & = \grd_{i,(\gra^{\grG^B_l})_0}\grd_{k,L_0(\grG^B_l)} \gra^{\grG^B_l} \\
\gra_{jh}^l & = \grd_{j,(\gra^{\grG^A_l})_1}\grd_{h,L_1(\grG^A_l)} \gra^{\grG^A_l}
\end{align*}
}
The equation \eqref{eq:puphigoF} follows now by the very definition.


\section{The action of the Weyl group}
For any $m\in P$ 
and for any $\grl \in Z$
we defined a variety $M_{m,\grl}(d,v)$. Observe that on both
$m,\grl$ there is a natural action of the Weyl group $W$. 
We define an action of the Weyl group also on $(d,v)$. We have already
described $d$ as an element of $X$ and $v$ as an element of $Q$. We can now
define
$$
\grs (d,v) = (d,\grs(v-d)+d).
$$
Observe that  $\grs(v-d)+d \in Q$ so the definition is well given.
So it make sense to consider the variety $M_{\grs m,\grs
\grl}(\grs(d,v))$ or the variety $\goM_{\grs \zeta} (  \grs(d,v))$
for $\zeta \in \goZ$.

In \cite{Na1} Nakajima used analytic methods to prove, in the
case of a finite Dynkin diagram, that if $\zeta$ is generic then there exists
a diffeomorphism of differentiable manifolds 
$$
\Phi_{\grs,\zeta} \colon  \goM_{\zeta} (d,v)
\lra
\goM_{\grs \zeta} (  \grs(d,v))
$$
and moreover that
$\Phi_{\grs',\grs\zeta}\comp\Phi_{\grs,\zeta}=\Phi_{\grs'\grs,\zeta}$.
In the same paper he also asserted that a similar construction could be obtained
in the general case using reflection functors as indeed we are going to do.

In \cite{Lu:Q4} Lusztig gave a purely algebraic construction of an isomorphism
$$
M_{0,\grl}(d,v) \isocan M_{0,s_i \grl} (s_i(d,v))
$$
whenever $\grl_i \neq 0$. In this paper we will give a generalization of Lusztig construction.

\begin{dfn}\label{dfn:22}
If $u \in \mZ^n=Q\cech$ and $A \subset Q\cech$ we define
$$
H_u = \{ (m,\grl)\in P \oplus Z \st \bra u, \cech ,\grl \ket = \bra u \cech , m \ket =0 \}
\;
\mand 
\;
H_A = \bigcup_{a\in A} H_a
$$
Let $K= \max \{1,a_{ij}^2 \st i,j \in I \}$. If $v \in \mZ^n$ we define
$$
\tilde U _{v} = \{ u \in \mN^I \tc 0 \leq u_i \leq K v_i \}
\; \mand \;
\tilde H^{v} = H_{\tilde U_v} .
$$
We define also
$$
 U_{\infty} = \bigcup _{i\in I} W \gra_i\cech \; \mand \; 
 H^{\infty} = H_{U_{\infty}}.
$$
Finally we set $ \calG_{v} = \{ (m,\grl) \in P \times Z_G \st \grs(m,\grl) \notin H^{\grs \cdot v} 
\text{ for all } \grs \in W \}.$
Both of the following definition of the set $\calG$ will be fine for us:
\begin{align*}
\calG &= \{(v,m,\grl)\in Q \times P \times Z_G  \st (m,\grl) \in \calG_{v}\} \; \text{ or } \\
\calG &= \{(v,m,\grl)\in Q \times P \times Z_G  \st (m,\grl) \notin H^{\infty} \}.
\end{align*}
\end{dfn}
We observe that in any case $\calG$ is $W$-stable.

\begin{prp}\label{prp:azioneweyl}
 For all $d,v$, for all $\grs \in W$  and for all $(m,\grl)$ such that 
$(m,\grl,v)  \in \calG_v$ there exists
an algebraic isomorphism:
$$
\Phi^{d,v}_{\grs,m,\grl}\colon M_{m,\grl}\bigl(d,v\bigr) 
\lra M_{\grs m ,\grs \grl }\bigl(\grs(d,v)\bigr).
$$
Moreover this isomorphims satisfies
\begin{equation}\label{eq:azioneweyl}
\Phi^{\grs(d,v)}_{\tau,\grs m,\grs \grl} \comp \Phi^{d,v}_{\grs,m,\grl} = 
\Phi^{d,v}_{\tau \grs,m,\grl}.
\end{equation}
\end{prp}


\subsection{Generators}
In this section we define the actions of the generators $s_i$ of $W$ following
\cite{Lu:Q4}. We fix $ i \in I$ and $(d,v)$, $\grl \in Z$ and $m \in P$. 
We call $(d,v')= s_i(d,v)$, $\grl ' = s_i \grl$
and $m '=s_i m$.
Through all this section we assume $v,v'\geq 0$.  
 For the convenience of 
the reader we write explicit formula in this case:
\begin{align*} 
   \grl_j ' &= \grl_j - c_{ij} \grl_i  &  m _j'&=  m_j - c_{ij} 
	 m_i  \; \text{ for all } j \\
   v_i' &= d_i -v_i +\sum _{j\neq i} a_{ij}  v_j && v_j'=v_j
	\; \text{ for all } j \neq i 
\end{align*}
Observe that we can choose 
$$ D_j' = D_j \;\text{ for all } j  \;\mand \;
	V_j' =V_j \;\text{ for all } j \neq i .
$$
In particular we have 
$$
  T_i = D_i \oplus \bigoplus _{h_1=i} V_{h_0} = T_i'
$$
since we suppose that  our quiver has not simple loops.

\begin{dfn}[Lusztig \cite{Lu:Q4}] \label{dfn:generatorsweyl}
Fix $\grl \in Z_G$ and define $Z^{\grl}_i (d,v)$
to be the subvariety of $S_i(d,v) \times S_i(d,v')$ of pairs $(s,s')= 
\bigl( (B,\grg,\grd), (B',\grg',\grd') \bigr)$ such that the following
conditions hold: 
\begin{enumerate}
   \item $B_h (s)= B_h'(s')$ for all $h$ such that $h_0,h_1 \neq i$, 
   \item $\grg_j (s)= \grg_j'(s') $ for all $j \neq i$ ,
   \item $\grd_j (s)= \grd_j'(s') $ for all $j \neq i$,
   \item the following sequence is exact:
\begin{equation} \label{eq:seqes:defsi}
\begin{CD} 
0 @>>> V_i' @>{a_i'}>> T_i @>{b_i}>> V_i @>>> 0,
\end{CD} 
\end{equation}
   \item $a_i'(s')b_i'(s') = a_i(s)b_i(s) - \grl_i \Id_{T_i}$,
   \item $s \in \grL_{\grl}(d,v)$ and $s' \in \grL_{\grl'}(d,v')$.
\end{enumerate}
\end{dfn}  

\begin{lem}
 Let   $(s,s') \in S_i(d,v) \times S_i(d,v')$ and suppose that it satisfies
conditions 1), 2), 3), 4) , 5) above then:
  \begin{enumerate}
	\item  $s \in \grL_{\grl}(d,v) \iff s' \in \grL_{\grl'}(d,v')$,
	\item if $\mu_j(s) - \grl_j \Id_{V_j}=0$ for all $j\neq i$ then $s \in \grL_{\grl}(d,v)$ ,
	\item if $\mu_j(s') = \grl_j \Id_{V_j'}$ for all $j\neq i$ then $s' \in \grL_{\grl'}(d,v')$.
   \end{enumerate} 
\end{lem}

\begin{proof}
2) We have to prove $b_ia_i - \grl_i\Id _{V_i}=0$ and by condition 4) it is enough to prove
 $b_ia_ib_i = \grl_i b_i$. So $b_ia_ib_i = b_i(a_i'b_i'-\grl_i) = \grl_i b_i$
by conditions 4) and 5).

The proof of 3) is equal to the proof of 2). We prove the implication
$\then$ in 1). By 2) and 3) it is enough to prove that $b_j'a_j' = \grl_j'$
for $j\neq i$. 
{\allowdisplaybreaks
\begin{align*}
 b_j'a_j' &= \sum_{h_1=j} \gre(h) B_h' B'_{\bar h} + \grg_j' \grd_j' = \\
	& = \sum_{h_1=j, h_0 \neq i} \gre(h) B_h B_{\bar h} + \grg_j \grd_j
		+ \sum_{h_1=j, h_0 = i} \gre(h) B_h' B_{\bar h}' = \\
	& = b_j a_j + \sum_{h_1=j, h_0 = i} \gre(h)
		\left( B_h' B_{\bar h}' - B_h B_{\bar h} \right) \\
	& = b_j a_j + \sum_{h_0=j, h_1 = i}  
		\left( B_{\bar h} \gre(h) B_{h}- B_{\bar h}' \gre(h) B'_h
		\right)  \\
	& = b_j a_j + \sum_{h_0=j, h_1 = i} \left( [a_ib_i]^{V_{h_0}}_{V_{h_0}}
  		- [a'_ib'_i]^{V_{h_0}}_{V_{h_0}} \right) \\
	& = \grl_j + \sum_{h_0=j, h_1 = i} \grl_i = \grl_j'
\end{align*}
}
The proof of the converse is completely analougous. 
\end{proof}

\begin{lem}
  Let $\grl \in Z_G$, $(s,s') \in Z^{\grl}_i(d,v)$ and  $\gra $ be 
  an element of the path algebra algebra of type $(\gra_0,\gra_1)$ then 
\begin{enumerate}
	\item if $\gra_0,\gra_1\neq i$ there exists an element $\gra'$ 
		of the $b$-path algebra of type  $(\gra_0,\gra_1)$ such that
	 	$\gra'(s')=\gra(s)$, 	 	
	\item if $\gra_1\neq i$ there exists an element $\gra'$ 
		of the $b$-path algebra of type  $(\gra_0,\gra_1)$ such that
		$\gra'(s') \grg '_{\gra_0} =\gra(s)\grg_{\gra_0}$, 	
	\item if $\gra_0\neq i$ there exists an element $\gra'$ 
		of the $b$-path algebra of type  $(\gra_0,\gra_1)$ such that
		$\grd '_{\gra_1}\gra'(s')  =\grd _{\gra_1}\gra(s)$, 	
	\item there exists an element $\gra'$ 
		of the $b$-path algebra of type  $(\gra_0,\gra_1)$ such that	
		$\grd '_{\gra_1}\gra'(s') \grg '_{\gra_0} = 	\grd
		_{\gra_1}\gra(s)\grg_{\gra_0}$, 
\end{enumerate} 
\end{lem}
\begin{proof}
By induction on the length of $\gra$ we can reduce the proof of this lemma 
to the following identities that are a consequence of condition 5) in
definition \ref{dfn:generatorsweyl}:
\begin{align*} B'_hB'_k &= 
\begin{cases}
  B_hB_k &\mif h \neq \bar k \\
  B_hB_k - \grl_i &\mif k = \bar h
\end{cases} \\
\grd '_iB'_k &=  \grd _iB_k \\
B'_h \grg '_i &= B_h \grg_i \\
\grd '_i \grg _i ' &= \grd _i \grg_i -\grl_i
\end{align*}
for $h,k$ such that $h_0 = i = k_1$.
\end{proof}

\begin{lem}
Let $(s,s') \in Z^{\grl}_i(d,v)$ and suppose $m_i \geq 0$ or $\grl_i\neq 0$ then
$$
s \text{ is } \chi_m \text{ semistable }
\iff
s' \text{ is } \chi_{m'} \text{ semistable }
$$
\end{lem}

\begin{proof}We prove only $\then$.
Let's do first the case $m_i \geq 0$.
If $s$ is $\chi_m$ semistable, then there exists $\grD= \{
A,B,\gra^{*}_{*} \}$ $m$-good  such that
$f_{\grD}(s) \neq 0$. 
Using the notation  in  \ref{sssec:generators}
we have $\grf_{\grD} = \det \Psi_{\grD}$ where
$ \Psi_{\grD} : Y \lra Z$ is a linear map. 
In our case we can write $Z$ as $\mC^{m_i} \otimes V_i
\oplus \wt Z $ and we obseve 
that no $V_i$
summunds appear in $Y$ or $\wt Z$. 

Now we
construct a new data $\grD'= \{ A',B',{\gra'}^*_* \}$ such that
$f_{\grD'}(s') \neq 0 $ and $f_{\grD'}$ a $\chi'$-covariant
polynomial.
Our strategy will be  the following: we substitute each $V_i$ with the space $T_i$
in the space $Z$ and we add $m_i$ copies of $V_i'$ to $Y$.
Let's do it more precise: 
first of all the new data will not be $m'$ good so  we have to
define  ${m'}_j^+  $ and ${m'}_j^-$:
\begin{enumerate}
	\item ${m'}_i^+ = 0$ and ${m'}_i^- = m_i = m_i^+$,
	\item ${m'_j}^- = m_j^-$ and ${m'_j}^+ = m_j^+ + a_{ij} m_i^+$ for all
		$j\neq i$, 	
	\item ${m'}^- = m^- $ and ${m'}^+ = m^+ + d_i m_i^+$.
\end{enumerate}
Observe that  ${m'_j}^+ - {m'_j}^- =  m _j'$ for all $j$ so
our data will furnish a $\chi'$ equivariant function.
Moreover if we define
$$
Z' = \mC^{m_i} \otimes V_i \oplus \wt Z  \;\mand\; 
Y' = \mC^{m_i} \otimes T_i \oplus Y 
$$
we observe that they have the numbers of $V'_j$, $\mC_a$, $\mC_b$ factors
specified by ${m'}$. Now we construct the new data $\grD'$ in such a way that 
with respect to the decompositions above we have:
\begin{align*}
 [\Psi_{\grD}(s)]_{\mC^{m_i}\otimes V_i \oplus \wt Z}^{Y} & =
\begin{pmatrix}
	( \Id\otimes b_i ) \comp \pi\\
	\Phi
\end{pmatrix},\\
 [\Psi_{\grD'}(s')]_{\mC^{m_i}\otimes T_i \oplus \wt Z}^{\mC^{m_i}\otimes
	V_i' \oplus Y} & = 
	\begin{pmatrix}
	\Id \otimes a'_i & \pi \\
	            0   &  \Phi
	\end{pmatrix}.
\end{align*}

If we construct a data with this property we observe that
$\Psi_{\grD}(s)$  is an isomorphism if and only if  
$\Psi_{\grD'}(s')$  is an isomorphism. Hence $f_{\grD}(s) \neq 0$ implies 
$f_{\grD'}(s') \neq 0$ and the lemma is proved.

To construct the new data we choose a basis 
$e_1,\dots,e_{d_i}$ of $D_i$ and we define the other elements of the data according to
the following rules 
\begin{enumerate}
	\item $A' =A$, 
	\item if  $B=(b_1,\dots,b_{m^+})$ we set $B'=(b_1,\dots,b_{m^+},
		\underbrace{e_1,\dots,e_1}_{m_1\;times},\dots ,
 		\underbrace{e_{d_i}, \dots e_{d_i}} _{m_1\;times})$, 
	\item ${\gra'}^{j_1,h_1}_{j_2,h_2}$ for $j_1 \neq i$ and 
		$h_2 \leq m^+_{j_2}$ is an
		element constructed according to case 1) in the
		previous lemma, 	
	\item ${\gra'}^{a}_{j_2,h_2}$ for $h_2 \leq m^+_{j_2}$  is an
		element constructed according to case 2) in the
		previous lemma, 
	\item ${\gra'}^{j_1,h_1}_{b_l}$ for $j_1 \neq i$ and 
		$l \leq m^+$ is an
		element constructed according to case 3) in the
		previous lemma.
	\item ${\gra'}^{i,h}_{b_l} = {\gra'}^{i,h}_{j,k}=0$
		if $l \leq m^+$ and $k \leq m_j^+$.
\end{enumerate} 
In this way we garantee that the projection of $\Psi_{\grD'} (s') $ onto 
$\wt Z$ is equal to $\bigl( 0\; \Phi  \bigr)$. To define the remaining  part
of the new data we do not give details on the indexes, but we explain how to
construct it.
It is  clear that we can choose ${\gra'}^{i,h}_*$ for the
remaining indeces * in such a way that the projection of 
$\Psi_{\grD'} (s')\bigr|_{ \mC^{m_i}\otimes V_i'}$ on
$\mC^{m_i} \otimes T_i$
is equal to $\Id \otimes a'_i$. 
Finally we
observe that a path $\grb$ from $V_j$ to $V_i$ with $j\neq i$ has to go
through a summand of $T_i$ so there exists a path $\gra$ such that 
$\grb(s)= b_i \comp \gra(s)$. Now we use the previous lemma to change $\gra$
with a $\gra'$  such that $\grb(s) = b_i \comp \gra'(s')$. 
More generally if $\grb$ is an element of the path algebra of type 
$(j,i)$ with $j\neq i$ then there exists an element
of the $b$-path algebra $\gra'$ such that 
$\grb(s) = b_i\comp \gra'(s)$.
In this way we
define the elements of the $b$-path algebra  
connecting summunds of $Y$ and summunds of $\mC^{m_i}\otimes T_i$.

In the case  $m_i<0$ we proceed in a similar way: we choose $\grD$ $m$-good
and we have 
$$
Y = \mC^{-m_i} \otimes V_i \oplus \wt Y, \quad
Y'=  \mC^{-m_i} \otimes T_i \oplus \wt Y, \quad
Z'= \mC^{-m_i} \otimes V'_i \oplus Z.
$$
As in the previous case we can find a new data $\grD'$ such that:
\begin{align*}
[\Psi_{\grD}(s)] _{ Z} ^{\mC^{-m_i}\otimes V_i \oplus \wt Y} &=
	\begin{pmatrix}
	\pi \comp  ( \Id \otimes a_i )  &
	\Phi
	\end{pmatrix},\\
 [\Psi_{\grD'}(s')] _{\mC^{-m_i}\otimes V_i' \oplus Z} ^{\mC^{-m_i}\otimes T_i
	\oplus \wt Y} & =  	
	\begin{pmatrix}
	\Id \otimes b'_i & 0 \\
	            \pi   &  \Phi
	\end{pmatrix}.
\end{align*}
Now to conclude that $\Psi_{\grD'}(s')$  is an isomorphism if $\Psi_{\grD}(s)$ is
we need to know that $b'_i$ is an
epimorphism and this is not garantee by $(s,s') \in Z^{\grl}_i(d,v)$. But if
$\grl_i \neq 0$ then, since $b_i'a_i'= -\grl_i$,  we have that
$b_i'$ is surjective.
\end{proof}

\begin{dfn}
 Let $p$ (resp. $p'$) be the projections of 
$Z_i^{\grl}(d,v)$ on $\grL_{\grl}(d,v) \subset S(d,v)$ 
(resp.  $\grL_{\grl'}(d,v') \subset S(d,v')$). Suppose that
$m_i > 0$ or $\grl_i \neq 0$ then we define  
$$
Z_i^{m,\grl}= p^{-1}\bigl( \grL_{m,\grl}(d,v) \bigr) 
= {p'}^{-1}\bigl( \grL_{m',\grl'}(d,v') \bigr).
$$
We define also 
$$
G_{i,v} = \prod _{j \neq i} GL(V_j) \times GL(V_i) \times GL(V_i').
$$
Observe that there are natural projections from 
$G_{i,v}$ to $G_v$ and $G_{v'}$, therefore ther are
natural actions of $G_{i,v}$ on $S_i(d,v)$, 
$S_i(d,v')$.
Observe that here is a natural action of 
$G_{i,v}$ on $Z_i^{\grl}$ and $Z_i^{m,\grl}$ such that 
the projections $p$, $p'$ are equivariant.
\end{dfn}

\begin{lem} \label{lem:abmonoepi}
  Let $s \in \grL_{\grl,m}(d,v)$ then 
\begin{enumerate}
	\item if $\grl_i \neq 0$ then $b_i$ is epi and $ a_i$ is mono,
	\item if $m_i>0$ then $b_i$ is epi,
	\item if $m_i<0$ then $a_i$ is mono.
\end{enumerate}
\end{lem}
\begin{proof}
   If $\grl_i \neq 0$ then the result  is clear by $b_ia_i = \grl_i$.
Suppose now that $\grl_i = 0 $ and $m_i >0$. Let $U_i = \Im b_i$ and let 
$V_i = U_i \oplus W_i$. Define now a one parameter subgroup $g(t)$ of  
$G_V$ in the following way: 
$$
[g_i(t)]^{U_i \oplus W_i}_{U_i\oplus W_i} = 
\begin{pmatrix}
1 & 0 \\ 0 &  t^{-1}
\end{pmatrix} 
\mand g_j \coinc 1 \mfor j
\neq i 
$$
Since $\Im b_i \subset U_i$ we have that there exists the limit 
$\lim_{t\to 0 } g(t) \cdot s = s_{0}$. Let now $n>0$ and $f$ a
$\chi^n$-covariant  function on $S$ such that $f(s)\neq 0$. Then
$$
f(s_0) = \lim_{t \to 0} f(g(t) \cdot s) = \lim_{t \to 0} 
{\det}_{GL(V_i)}^{nm_i} f(s) = \lim_{t \to 0} 
t^{-nm_i\dim W_i} f(s)
$$
So we must have $\dim W_i =0$.
The proof of the third case is completely similar to this one.
\end{proof}

\begin{lem}[see also Lusztig \cite{Lu:Q4}]

If $m_i>0$ or $\grl_i \neq 0 $ then 
\begin{enumerate}
	\item $p: Z_i^{m,\grl}(d,v) \lra \grL_{m,\grl}(d,v)$ is a principal 
		$GL(V_i')$ bundle,
	\item $p': Z_i^{m,\grl}(d,v) \lra \grL_{m',\grl'}(d,v')$ is a principal 
		$GL(V_i)$ bundle.
\end{enumerate} 
\end{lem}
\begin{proof}
Lusztig's proof extend to this case without changes. 
Let's prove for example 1. We have to prove: $i.$ that the action on the fiber is free,
$ii.$ that it is transitive. First of all we observe that by the previous lemma 
if $s \in \grL_{m,\grl}$ then $b_i(s)$ is epi. In particular there exists
$a_i' : V_i' \lra T_i$ such that sequence \eqref{eq:seqes:defsi} is exact, and clearly $a_i'$ 
is univoquely determined up to the action of $GL(V_i')$, moreover this action is free.
So $i.$ and $ii.$ reduce to the following fact: if $s \in \grL_{m,\grl}$ and $a_i'$
is such that  sequence \eqref{eq:seqes:defsi} is exact, then there exists 
a unique, $b_i'$ such that $a_i'b_i'=a_ib_i -\grl_i$. Since $a_i'$ is mono the unicity is clear.
To prove the existence we observe that it is equivalent to 
$\Im a_i' \supset \Im (a_ib_i - \grl_i)$. But the last statement is clear since we have:
$\Im a_i' = \ker b_i$ and $ b_i(a_ib_i - \grl_i) = 0$.
\end{proof}

\begin{prp} \label{prp:isoMsi}
 If $m_i>0$ or $\grl_i \neq 0$ then the projections $p$, $p'$ induces 
algebraic isomorphisms $\bar p$ , $\bar p'$:
$$
\begin{CD}
\grL_{m,\grl}(d,v) /\!/ G_v 
@<{\circa}<{\bar p}< 
Z_i^{m,\grl}(d,v) /\!/ G_{i,v} 
@>{\circa}>{\bar p '}> 
\grL_{m',\grl'}(d,v') /\!/ G_{v'}
\end{CD}
$$
\end{prp}
\begin{proof}
This proposition is a straightforward consequence of the previous lemma and the following general
fact (see for example \cite{GIT} Proposition 0.2): 
let $G$ be an algebraic groups over $\mC$ and  $X$, $Y$ two irreducible 
algebraic variety over $\mC$; if $G$ acts on $X$ and $\grf:X \lra Y$ is 
such that for all $y \in Y$ the fiber $X_y$ contains exactly one $G$-orbit
then $\grf $ is a categorical quotient.
If we apply this lemma to the projection $p$, (resp. $p'$)  and to the group 
$GL(V_i')$ (resp. $GL(V_i)$) we obtain the required result.
\end{proof}

We can use this proposition to define the action of the generators of the Weyl group.
\begin{dfn}
Let $i,\grl,m,d,v,\grl',m',v'$ be as above, and suppose 
$d_j \geq 0$, $ v_j,v_j' \geq 0 $ for all $j$ then we define an isomorphism 
of algebraic variety
$$
\Phi^{d,v}_{s_i,\grl,m}\colon M_{m,\grl}(d,v) \lra M_{m',\grl'}(d,v')
$$
in the following way:
\begin{enumerate}
	\item if $m_i>0$ or $\grl_i \neq 0$ 
		then we set $\Phi^{d,v}_{s_i,\grl,m} = \bar p' {\bar p}^{-1}$,
	\item if $m_i<0$ then we exchange the role of $v,v'$ in the previous construction: more
	precisely we observe that $m_i' >0$ so we can define 
	$\Phi^{d,v'}_{s_i,\grl',m'}\colon M_{m',\grl'}(d,v') \lra M_{m,\grl}(d,v)$
	and we define $\Phi^{d,v}_{s_i,\grl,m} = \bigl(\Phi^{d,v'}_{s_i,\grl',m'}\bigr)^{-1}$.
\end{enumerate}
\end{dfn}

\begin{oss}\label{oss:defisosi}
To see that $\Phi^{d,v}_{s_i,\grl,m}$ is univoquely defined we have to verify 
that if $\grl_i \neq 0$ and $m_i<0$ the two definitions above coincide. This fact reduces easily
to the following remark: if $\grl_i \neq 0$ then 
$$
(s,s') \in Z_i^{\grl}(d,v) \iff (s',s) \in Z_i^{\grl'}(d,v').
$$
Let us prove, for example, the $\then$ part. 
Since $a_i b_i = a_i'b_i' + \grl_i = a_i'b_i' - \grl_i ' $ the only thing we have to verify
is that the sequence
$$
\begin{CD}
0 @>>> V_i @>{a_i}>> T_i @>{b_i'}>> V_i' @>>> 0
\end{CD}
$$
is exact. The surjectivity of $b_i'$ and the injectivity of $a_i$ are a consequence of
$\grl_i \neq 0$. Since $\dim T_i = \dim V_i + \dim V_i'$ we need only to prove that 
$b_i'a_i =0$. Observe that $b_i'a_i =0$ if and only if $ a_i'b_i'a_i =0$ since also $a_i'$ 
is injective. Finally $ a_i'b_i'a_i = (a_ib_i-\grl_i) a_i = 0$.
\end{oss}


\subsection{Preliminaries} We saw how to define 
$$
\Phi^{d,v}_{s_i,m,\grl}\colon M_{m,\grl} \bigl(d,v\bigr)\lra 
M_{s_i(m),s_i(\grl)}\bigl( s_i(d,v) \bigr)
$$
in the case that $(\grl_i,m_i ) \neq 0$ and $d,v,s_iv \geq 0$. 
To define an action of the Weyl group
we have now to garantee that coxeter relations hold. We will prove these relations 
in the next paragraph. Before doing it we observe that we have to garantee some conditions 
on $m,\grl$ such that we will be able to define 
$\Phi_{s_i,\grs m,\grs\grl}^{\grs(d,v)}$ for any element
$\grs \in W$: this condition will be $(m,\grl) \in \calG_v$ (\ref{dfn:22}).
We have also to say something about the case $d_i <0$ or $v_i<0$ for some 
$i \in I$. 

In the case that $d_i <0$ for some $i$ then $M_{m,\grl}\bigl( \grs(d,v) \bigr)= \vuoto
$ for all $\grs,m,\grl$ by the very definition, so there is nothing to define.

The second trivial case is 
$d=v=0$. Indeed in 
this case we have $M_{\grs m,\grs \grl}\bigl(\grs(d,v)\bigr) = \{ 0 \}$
so the definition is trivial. 

The other two cases are threated in the two lemmas below.

In the following we fix $d$ such that $d_i \geq 0 $ for all $i$. It will be convenient to 
define an affine action of $W$ on $Q$ by $\grs \cdot v = \grs(v-d) +d$.

\begin{lem}
Let $d\geq 0 $ and $(m,\grl) \in \calG_v$  if there exists $\grs$ such that 
$\grs \cdot v \not \geq 0$ then
$M_{m,\grl}(d,v) = \vuoto$.
\end{lem}
\begin{proof}
Suppose that $\grs$ is an element of minimal length such that 
$\grs \cdot v \not \geq 0$ and let $l = \ell(\grs)$. We prove the lemma by induction
on $l$. The case $l=0$ is trivial.

\noindent \emph{Initial step:  $l=1$.} If $s_i \cdot v \not \geq 0$ then we have
$0 \leq d_i + \sum a_{ij} v_j < v_i$. Hence $\dim T_i < \dim V_i $,
$u = ( 0,\dots,\overset{i}{1},0,\dots) \in \tilde U_v$ and $(\grl_i,m_i) \neq 0$. 
So $M_{m,\grl}(d,v) = \vuoto$ by lemma \ref{lem:abmonoepi}.

\noindent \emph{Inductive step: if $l\geq2$ then $l-1 \then l$.} Let $\grs = \tau s_i$ with
$\ell(\tau) = l-1$ and 
$v' = s_i \cdot v$, $\grl'= s_i \grl$, $m ' = s_i m$. By induction
$M_{m',\grl'}(d,v') = \vuoto$ and, since $l \geq 2 $, $v' \geq 0$.
If $(m_i,\grl_i) \neq 0$ then we can apply Proposition 
\ref{prp:isoMsi} and we obtain $M_{m,\grl}(d,v)\isocan
M_{m',\grl'}(d,v') = \vuoto$. If $(m_i,\grl_i)=0$ then
$u = ( 0,\dots,\overset{i}{1},0,\dots) \notin \tilde U_v$, hence
$v_i =0$. Moreover $\grl ' = \grl$ and $m' = m$ so 
$(m_i',\grl_i')=0$ and $u = ( 0,\dots,\overset{i}{1},0,\dots) 
\notin U_{v'}$. Hence $v_i'=0$ so $v'=v$ and $\tau v \not \geq 0$ against the minimality
of $\grs$.
\end{proof}

\begin{lem}\label{lem:sisipuodefinire}
 Let $(I,H)$ be connected, $(m,\grl ) \in \calG_v$ and suppose $d\geq 0$ and 
$\grs \cdot v \geq 0$ for all $\grs \in W$. 
If there exists $i \in I, \grs \in W$ such that $\grs(m,\grl) = ( m',\grl')$ and 
$(m_i',\grl_i')=0$ then $d=v=0$.
\end{lem}  
\begin{proof}
Without loss of generality we can assume $\grs =1$. 

\noindent \emph{First step}: $v_i =0$. This is clear since otherwise $u = ( 0,\dots,\overset{i}{1},0,\dots) \in U_v$.

\noindent \emph{Second step}: $d_i =0$ and $v_j=0$ for all $j$ such that $a_{ij}$. 
Let $v'= s_i\cdot v$ and observe that $s_i\grl =\grl$ and $s_im=m$. Then as in first step 
we have $0 =v_i' = d_i + \sum_{j}a_{ij}v_j$ from which the claim follows.

Let now $ W' = \bra \{s_j \st a_{ij} \neq 0 \mand j\neq i\} \ket $. 
If $(d,v) \neq 0$ then there exists 
$j \in I$ and $\grs \in W'$ such that $a_{ij} \neq 0$ and
$$
n = d_j + \sum_{h\in I} a_{jh}\Tilde v _h >0.
$$
where $\Tilde v=\grs\cdot v$. Since $(\grs \grl)_i = \grl_i = 0 = m_i = (\grs m)_i$ we can assume
$\grs = 1$. Let now 
$v' = s_is_j \cdot v$, $ \grl' =s_is_j \grl$ and $m ' = s_is_j m$, we have:
\begin{align*}
	v'_i &= a_{ij}n & \grl'_i&= -a_{ij}\grl_j &  m_i&= -a_{ij}m_j \\
	v'_j &= n & \grl'_j&= (a_{ij}^2-1) \grl_j &  m_j&= (a_{ij}^2-1)m_j .
\end{align*}
Hence  $u= (0,\dots,\overset{j}{a_{ij}},0,\dots, \overset{i}{a_{ij}^2-1},0,\dots)\in U_{v'}$ and
$\bra u \cech ,\grl' \ket = \bra u\cech, m' \ket =0$ against $(m,\grl) \in \calG_v$.
\end{proof}

\begin{oss}
the analougous lemma in the case 
of $\calG= \{(m,\grl,v) \st (m,grl) \notin H^{\infty} \}$ are more simple. 
\end{oss}


\subsection{Relations}
In this section we define an isomorphism of algebraic variety
$$
\Phi_{\grs,m,\grl}^{d,v} \colon 
M_{m,\grl}\bigl(d,v\bigr) \lra M_{\grs m,\grs\grl}\bigl(\grs(d,v)\bigr).
$$
in the case $(m,\grl) \in \calG _v$ or $(m,\grl) \notin H^{\infty}$.
In the case $d \not \geq 0$ or in the case in which there exists $\grs \in W$
such that  $\grs v \not \geq 0$ or in the case $d=v=0$ we have seen in the previous  
section
that there is nothing to define or that the definition is trivial. 
In the remaing cases we observe that 
for all $\tau, i$ we have $(\tau(m)_i , \tau(\grl)_i) \neq 0$ by lemma
\ref{lem:sisipuodefinire}. Hence we can define $\Phi_{\grs,m,\grl}^{d,v}$
by induction on $\ell(\grs)$ by the formula
\begin{equation} \label{eq:defphisigma}
\Phi_{\grs,m,\grl}^{d,v} = 
\Phi_{s_i,\tau m,\tau \grl}^{\tau(d,v)}
\comp
\Phi_{\tau , m,\grl}^{d,v}.
\end{equation}
Of course we have to prove that this definition is well given by checking
Coxeter relations:
$$
	s_i^2 = \Id, \quad
	s_is_j=s_js_i \;\mif a_{ij}=0  \; \mand \;
	s_is_js_i = s_js_is_j\;\mif a_{ij}=1 
$$
which in our situation take the following form:
\begin{subequations}
\begin{gather}
\Phi_{s_i,s_i\grl,s_im}^{s_i(d,v)} \comp \Phi_{s_i,\grl,m}^{d,v} = \Id \label{eq:cox1}\\
\Phi_{s_i,s_j\grl,s_jm}^{s_j(d,v)} \comp \Phi_{s_j,\grl,m}^{d,v} =
\Phi_{s_j,s_i\grl,s_im}^{s_i(d,v)} \comp \Phi_{s_i,\grl,m}^{d,v} \label{eq:cox2}\\
\Phi_{s_i,s_js_i\grl,s_js_im}^{s_js_i(d,v)} \comp 
\Phi_{s_j,s_i\grl,s_im}^{s_i(d,v)} \comp
\Phi_{s_i,\grl,m}^{d,v}
=
\Phi_{s_j,s_is_j\grl,s_is_jm}^{s_is_j(d,v)} \comp
\Phi_{s_i,s_j\grl,s_jm}^{s_j(d,v)} \comp
\Phi_{s_j,\grl,m}^{d,v}.  \label{eq:cox3}
\end{gather}
\end{subequations}
The first of the two equations is clear by the very definition and remark
\ref{oss:defisosi}. The second equation is trivial. We need to prove the third equation.
We will need the following two simple lemmas of linear algebra which proofs are trivial.

\begin{lem}\label{lem:25}
Let $V,W,X,Y,Z$ be finite dimensional vector spaces and $\gra,\grb,\grg,\grd,\gre,\grf$
linear maps between them as in the diagrams below.
The diagram 
$$
\begin{CD}
0 @>>> V @>{\begin{pmatrix} \gra \\ \grb \\ \grg \end{pmatrix}}>> 
W \oplus X \oplus Y @>{\begin{pmatrix} \grd & 0 & -1 \\ 0 & \gre & \grf \end{pmatrix}}>>
Y \oplus Z @>>> 0
\end{CD}
$$ 
is exact if and only if the diagram 
$$
\begin{CD}
0 @>>> V @>{\begin{pmatrix} \gra \\ \grb \end{pmatrix}}>> 
W \oplus X @>{\begin{pmatrix} \grf \grd  & \gre \end{pmatrix}}>>
Z @>>> 0
\end{CD}
$$ 
is exact and $\grg = \grd\gra$.
\end{lem}

\begin{lem}\label{lem:26}
Let $U,V,W,X,Y,Z$ be finite dimensional vector spaces and 
$\gra,\grb,\grg,\grd,\gre,\grf,\psi,\rho,\grs$
linear maps between them as in the diagrams below such that
$\psi \oplus \rho : W \oplus X \lra Z$ is an
epimorphism.
Then the  diagram 
$$
\begin{CD}
0 \lra U 
@>{\begin{pmatrix} \gra \\ \grb \\ \grg \end{pmatrix}}>> 
V \oplus W \oplus X 
@>{\begin{pmatrix}  \grd & 0 & 1 \\ \gre & \grf & 0 \\ 0 & \psi & \rho \end{pmatrix}}>>
X \oplus Y \oplus Z
@>{\begin{pmatrix}\rho & \grs & -1 \end{pmatrix}}>>
Z \lra 0
\end{CD}
$$ 
is exact if and only if $\grg = - \grd \gra$ , $\psi=\grs\phi$ , 
$\rho\grd + \grs\gre =0$ and the diagram
$$
\begin{CD}
0 @>>> U @>{\begin{pmatrix} \gra \\ \grb \end{pmatrix}}>> 
V \oplus W @>{\begin{pmatrix} \gre & \grf \end{pmatrix}}>>
Y @>>> 0
\end{CD}
$$
is exact.
\end{lem}

We fix now and  $i,j$ such that $a_{ij}=1$ and we verifies \eqref{eq:cox3}. 
Let 
\begin{align*}
 \grl '  & = s_i \grl    & m'    &=s_i m     & v'    & =s_i v     \\
 \grlsec & = s_j \grl '  & \msec &=s_j m'    & \vsec & =s_j v'    \\ 
 \grlter & = s_i \grlsec & \mter &=s_i \msec & \vter &=s_i \vsec  \\
 \Tilde {\grl}  & = s_j \grl  & \Tilde m   &=s_j m & \Tilde v  & = s_j v \\
 {\Tilde {\Tilde {\grl}}} & = s_i {\Tilde {\grl}} 
	& {\Tilde {\Tilde m}} &=s_i {\Tilde m} & {\Tilde {\Tilde v}} & = s_i {\Tilde v }
\end{align*}

First of all we observe that since relation \eqref{eq:cox1} holds we can assume that:
\begin{enumerate}
	\item $\grl_i \neq 0$ or $m_i > 0$ and $\grl_j \neq 0 $ or $m_j >0$, 
	\item $\grl'_j \neq 0$ or $m'_j > 0$ and $\Tilde{\grl} _i \neq 0 $ or $\Tilde m _i >0$, 
	\item $\grlsec_i \neq 0$ or $\msec_i > 0$ and 
		$\Tilde{ \Tilde{\grl}}_j \neq 0 $ or $\Tilde{\Tilde m}_j >0$. 
\end{enumerate}
Define 
\begin{align*}
Z_{iji} & = \{ (\ster,s) \in \grL_{\mter,\grlter} (d,\vter) \times \grL_{m,\grl}(d,v) \st
\exists \ssec \in S(d,\vsec),  \\ 
& \;\qquad \mand s' \in S(d,v') \msuchthat   
(\ster,\ssec) \in Z_i^{\msec,\grlsec}(d,\vsec), \\
& \;\qquad (\ssec,s')\in Z_j^{m',\grl'}(d,v') \mand
(s',s) \in Z_i^{m,\grl}(d,v) \} \\
Z_{jij} & = \{ (\ster,s) \in \grL_{\mter,\grlter} (d,\vter) \times \grL_{m,\grl}(d,v) \st
\exists \Tilde{\Tilde{s}} \in S(d,\Tilde{\Tilde{v}}),  \\ 
& \;\qquad \mand \Tilde s \in S(d,\Tilde v ) \msuchthat   
(\ster,\Tilde{\Tilde s}) \in Z_j^{\Tilde{\Tilde m},\Tilde{\Tilde{\grl}  } }(d,\Tilde{\Tilde v}), \\
& \;\qquad (\Tilde{\Tilde s},\Tilde s )\in Z_i^{\Tilde m,\Tilde{\grl}}(d,\Tilde v ) \mand
(\Tilde s,s) \in Z_j^{m,\grl}(d,v) \} 
\end{align*}
Observe that $(\ster,s) \in Z_{iji} \iff p^{d,\vter}_{\mter.\grlter}(\ster) = 
\Phi_{s_i}\Phi_{s_j}\Phi_{s_i}(p^{d,v}_{m,\grl}(s))$ and that 
$(\ster,s) \in Z_{jij} \iff p^{d,\vter}_{\mter.\grlter}(\ster) = 
\Phi_{s_j}\Phi_{s_i}\Phi_{s_j}(p^{d,v}_{m,\grl}(s))$. So 
relation \eqref{eq:cox3} is equivalent to $Z_{iji}=Z_{jij}$.

Let now $R_i = D_i \oplus \bigoplus_{h\st h_1=i, h_0 \neq j} V_{h_0} $, 
$R_j = D_i \oplus \bigoplus_{h\st h_1=i, h_0 \neq j} V_{h_0} $ and observe 
that $T_i = R_i \oplus V_j$ and $T_j = R_j \oplus V_i$.  Let 
$k$ be the only element of $H$ such that $k_0=j$ and $k_1 = i$. Let $\gre = 
\gre(k)$. Define also $A=A(s)= B_k(s)$, $B=B(s)=B_{\bar k}(s)$ and for 
$l =i ,j $ and $\{l',l\}=\{i,j\}$ set 
$c_l=c_l(s)= \pi_{R_l}^{R_l\oplus V_{l'}} a_l(s)$ and  
$d_l=d_l(s)=  b_l(s)\bigr|_{R_l} $. 

Let now $(s,s\ter) \in \grL_{\grl}(d,v) \times \grL_{\grl\ter}(d,v\ter)$ and set 
$A^* =A(s^*)$,
$ B^* =B(s^*)$,
$ c_l^* = c_l(s^*)$ and
$ d_l^* = d_l(s^*)$ for $l \in \{i,j\}$ and $* \in \{\;, \ter \}$.

If we apply lemmas \ref{lem:25} and \ref{lem:26} to our situation we obtain
the following result: $(s,\ster) \in Z_{iji} $ of and only there exists 
vector spaces $V_i',V_j',V_i\sec ,V_j\sec $ and linear maps
$A',B',c_i',d_i',c_j',d_j', 
A\sec , B\sec , c_i\sec , d_i\sec , c_j\sec , d _j\sec$
sucht that:
\begin{enumerate}
 \item $\dim V_l^* =v_l^*$ for $l \in \{i,j\} $ and $* \in \{ ',\sec\}$,
 \item for each $*\in \{\, ' , \,\sec\}$ and 
	$l \in \{ i, j\}$
	$A^* \in \Hom (V^*_i,V^*_j)$, $B^* \in \Hom(V^*_j , V^* _i)$,
	$c_l \in \Hom (V^*_l,R^*_l)$ and 
	$d_l \in \Hom (R^*_l,V^*_l)$,
 \item $V_j \ter = V_j \sec$, $c_j \ter = c_j \sec$, $d_j \ter = d_j \sec$ 
	and 
	\begin{align*}
	c_i\ter d_i \ter &= c_i d_i -\grl_i -\grl_j &  c_i\ter B\ter& =c_i'B\sec \\
	A\ter d_i\ter &=A\sec d_i'& \gre A\ter B\ter&= \gre A\sec B\sec - \grl_j 
	\end{align*}
 \item $V_i \sec = V_i '$, $c_i \sec = c_i '$, $d_i \sec = d_i'$ 
	and 
	\begin{align*}
	c_j\sec d_j \sec &= c_j d_j -\grl_i -\grl_j &  c_j\sec A\sec & =c_j A' \\
	B\sec d_j\sec &= B' d_j & \gre A \sec B \sec&= \gre A' B' + \grl_i + \grl_j 
	\end{align*}
 \item $V_j' = V_j $, $c_j' = c_j$, $d_j ' = d_j$ 
	and 
	\begin{align*}
	c_i' d_i ' &= c_i d_i -\grl_i  &  c_i' B'& =c_i B \\
	A' d_i' &=A d_i & \gre A' B'&= \gre A B - \grl_i 
	\end{align*} 
 \item $ \gre c_i' B\sec A\ter + c_i' d_i' c_i \ter =0$ and $\gre A' B\sec = d_ja_j\sec$,
 \item the following diagrams are exact
$$
\begin{CD}
0@>>> V_i\ter
@>{\begin{pmatrix} c_i \ter \\ c_j\sec A\ter \end{pmatrix}}>>
R_i \oplus R_j
@>{\begin{pmatrix} Ad_i& d_j \end{pmatrix}}>> 
V_j @>>> 0 \\
0@>>> V_j\sec
@>{\begin{pmatrix} c_j \sec \\  c_i'B\sec \end{pmatrix}}>> 
R_j \oplus R_i
@>{\begin{pmatrix} Bd_j& d_i  \end{pmatrix}}>> 
V_i @>>> 0 \\
0@>>> V_i'
@>{\begin{pmatrix} c_i' \\  A' \end{pmatrix}}>> 
R_i \oplus V_j
@>{\begin{pmatrix} d_i & \gre B  \end{pmatrix}}>>
V_i @>>> 0 
\end{CD}
$$
\end{enumerate}
\begin{oss}
The first condition in point 6) is equivalent to $\gre B\sec A\ter = d_i' c_i \ter$. 
Indeed this condition is certainly sufficient. 
To prove the necessity observe that by the injectivity of $a_i' = (c_i '\;A')^t$ it is enough to 
prove $ \gre c_i' B\sec A\ter + c_i' d_i' c_i \ter =0$ and 
$ \gre A' B\sec A\ter + A' d_i' c_i \ter =0$. The first equation is the first condition
in point 6) and the second one is a consequence of 
$ \gre A'B\sec =d_j c_i \ter $, $A'd_i'=Ad_i$ and the exactness of the first sequence.
\end{oss}
\begin{oss}
The condition $(s,s\ter) \in Z_{jij}$ can be expressed in a similar way.
In the prevoius conditions we have only to change $i$ with $j$ and $\gre$ with 
$-\gre$.
\end{oss}

We will prove now $Z_{iji} \subset Z_{jij}$.
To do it we supose that 
$A', \dots, d_j \sec$ are given as above and we  construct
$\Tilde A ,\Tilde B, \Tilde c_i , \Tilde d_i ,\Tilde c_j, \Tilde d_j, 
\Tilde{\Tilde {A}} , \Tilde{\Tilde { B}} ,\Tilde{\Tilde { c}}_i ,\Tilde{\Tilde { d}}_i , 
\Tilde{\Tilde { c}}_j ,\Tilde{\Tilde { d}} _j$ such that they satisfy the conditions.
for $(s,s\ter) \in Z_{jij}$.

\noindent {First step: construction of 
$\Tilde A, \Tilde B, \Tilde c_i, \Tilde c_j , \Tilde d_i, \Tilde d_j$}. 
Choose $\Tilde s$ 
such that $(\Tilde s, s ) \in Z_j^{\chi,\grl}$
and define $\Tilde A = A(\Tilde s)$, 
$\Tilde B = B( \Tilde s)$, $\Tilde c_l = c_l(\Tilde s)$ and $\Tilde d_l = d_l(\Tilde s)$ for 
$l \in \{i,j\}$ .

Now I claim that there exists unique $\Tilde {\Tilde A} : V_i\ter \lra \Tilde V_j$ and 
$\Tilde {\Tilde B} : \Tilde V_j \lra  V_i\ter$ such that:
$$
\begin{cases}
 \Tilde c_j \Tilde{\Tilde A} = c_j\sec A\ter \\
 \Tilde B \Tilde{\Tilde A}   = - \gre d_i c_i \ter
\end{cases}
\;\mand\;
\begin{cases}
 \Tilde{\Tilde A} \Tilde{\Tilde B} = \Tilde A \Tilde B - \gre \grl_i - \gre \grl_j \\
 c_i \ter  \Tilde{\Tilde B} = c_i \Tilde B
\end{cases}
$$

\noindent \emph{Unicity of $\Tilde{\Tilde A}$}: since the map 
$\Tilde a_j = (\Tilde c_j \;\; - \gre \Tilde B )^t$
is injective the unicity is clear. 

\noindent \emph{Existence of $\Tilde{\Tilde A}$}: to prove the existence of $\Tilde{\Tilde A}$
is enough to prove:
$$
\Im 
\begin{pmatrix} c_j \sec  A \ter \\ - \gre d_i c_i \ter \end{pmatrix}
\subset 
\Im 
\begin{pmatrix} \Tilde c_j \\ \Tilde B \end{pmatrix} = 
\ker 
\begin{pmatrix} d_j & -\gre A \end{pmatrix}.
$$
So the thesis follows from $d_jc_j\sec A\ter + A d_i c_i \ter =0$.

Let now $\Tilde {\Tilde a} _i = (c_i \ter \;\; \Tilde {\Tilde A })^t$. I claim 
that  $\Tilde {\Tilde a} _i $ is injective and that $\Im  \Tilde {\Tilde a} _i = 
\ker ( d_i \;\; \gre \Tilde B)= \ker \Tilde b _i$. First of all observe that since 
$\Tilde m_i >0$ or $ \grl_i \neq 0 $, $\Tilde b_i$ is surjective. Observe also that 
$$
\begin{pmatrix}
\Tilde c_j & 0 \\ 0 & \Id_{V_i\ter}
\end{pmatrix}
\comp
\begin{pmatrix}
\Tilde {\Tilde A}  \\ c_i \ter  
\end{pmatrix}
=
\begin{pmatrix}
c_j \sec A\ter \\ c_i \ter  
\end{pmatrix}.
$$
So $\Tilde {\Tilde a}_i$ is injective as claimed.
Now since $\dim R_i + \dim \Tilde V_j = \dim V_i \ter + \dim V_i$ to prove the last part of the claim it is enough to check that $\Tilde b_i \Tilde{\Tilde a}_i=0$. Indeed 
$$
\Tilde b_i \Tilde{\Tilde a}_i = d_i c_i \ter +\gre \Tilde B \Tilde{\Tilde A} =0.
$$
\noindent \emph{Unicity of $\Tilde{\Tilde B}$}: this is a consequence of $\Tilde {\Tilde a}_i$
injective.

\noindent \emph{Existence of $\Tilde{\Tilde B}$}: As for the existence of $\Tilde{\Tilde A}$ 
this is equivalent to 
$$
\Im 
\begin{pmatrix} c_i \Tilde B \\ \Tilde A \Tilde B -\gre \grl_i - \gre \grl_j \end{pmatrix}
\subset 
\Im 
\begin{pmatrix} c_i \ter \\ \Tilde{\Tilde A}   \end{pmatrix} = 
\ker 
\begin{pmatrix} d_i & \gre \Tilde B \end{pmatrix}.
$$
So the thesis follows from
$ \gre \Tilde B \Tilde A \Tilde B - \grl_i \Tilde B - \grl_j \Tilde B + d_i c_i \Tilde B =0$.

Finally we set 
\begin{align*}
 \Tilde{\Tilde V}_i & = V_i \ter   & \Tilde{\Tilde c}_i &=c_i \ter    
& \Tilde{\Tilde d}_i&=d_i\ter \\
 \Tilde{\Tilde V}_i & =\Tilde V_j  & \Tilde{\Tilde c}_j &=\Tilde c_j  
& \Tilde{\Tilde d}_j&=\Tilde d_j.
\end{align*}
The verification of all the conditions is now straightforward.

The inclusion $Z_{jij}\subset Z_{iji}$ can be proved similarly and equation 
\eqref{eq:azioneweyl} is clear by definition. So Proposition \ref{prp:azioneweyl} is proved.

\section{A representation of the Weyl group}
In this section, following Nakajima \cite{Na1}, we show how to use the above action to
construct an action of the Weyl group on the homology of quiver varieties.
Maybe this action is related with the one constructed by Slodowy in the case of flag varieties 
(\cite{Slodowy}, ch.4).

First we recall some general about the action of the Weyl group. 
Let $Z\cech = Q\cech \otimes _{\mZ} \mC$ and $Z= Q \otimes _{\mZ}P$. On $Z$, 
$Z\cech$ there is a natural action of $W$.
\begin{lem}
For all $ u\in Z\cech$ the set $Wu$ is discrete.
\end{lem}
\begin{lem}
Consider the action of $W$ on $\mP(Z\cech)$. If $p \in \mP(Z\cech)$ then 
$$
\overline{Wp} \text{ is countable}.
$$
\end{lem}
If $p \in \mP(Z\cech)$ we define $H_p= \{x \in P\otimes _{\mZ}\mC \st \bra x , p \ket =0\}$
\begin{lem}
If $p \in \mP(Z\cech)$ then
$$
\overline{W H_p} = \bigcup_{q \in \overline{Wp}} H_q.
$$
\end{lem}

We define 
$\calH = \overline{W\calH_v \cup \calH_U}$ and  
$\calR = \goz - \calH$. 
By the previous lemmas $\calH$ is the union of a countable number of real codimension
$3$ subspaces in $\goZ$ and in particular $\calR$ is simply connected.
We need also the following definition
\begin{align*}
K &= \{ u \in \mZ^I  \st -\sum{u_i\gra_i} \text{ is dominant and supp}\, u \;
     \text{ is connected } \} \\
P_0 &= \{ p \in P \st p \text{ is dominant and } \bra u \cech , p \ket \geq 2 
	\text{ for all } u \in K \}.
\end{align*}

Now we choose $d,v$ such that $\bar d =\bar v$.

\begin{lem}\label{lem:lemma38}
If $\bar d \in P_0$ then $\tilde {\grm}$ 
is surjective and is a locally trivial bundle over
$\calR$.
\end{lem} 
\begin{proof}
By Proposition 10.5 and Corollary 10.6 in \cite{Na2} there exists a closed orbit $Gs$ in 
$\grL_0(d,v)$ with trivial stabilizer. Then by Proposition \ref{prp:KNM1} there exists 
$t \in Gs$ such that $\tilde {\grm}(t) = 0$ and by lemma 
\ref{lem:regmu} and \ref{regolaritamu}
$d\tilde {\mu}_t$ is surjective. 
Now the surjectivity follows by homogeneity.

The local triviality over $\calR$ follows also from lemma 
\ref{lem:regmu} and \ref{regolaritamu}.
\end{proof}

Now consider 
\begin{align*}
R & = \{\grl \in Z \st (0,\grl) \notin \calR \}, \\
\grL(d,v) & = \{ (\grl,s) \in Z \times S \st s \in \grL_{\grl} \} , \\
M(d,v) &= \grL(d,v) /\!/ G_v \;\mand\; p : M(d,v) \lra Z \text{ the projection} \\
\goL(d,v) & = \{ (\zeta,s) \in \goZ \times S \st s \in \goL_{\zeta} \} , \\
\goM(d,v) &= \goL(d,v) / U(V) \;\mand\; \tilde p : \goM(d,v) \lra \goZ \text{ the projection}
\end{align*}
We have the following commutative diagram
$$
\xymatrix{
(\grl,s) \ar@{}[r]|{\in} \ar[d] & M(d,v) \ar[r]^p \ar[d]^{\mi_M} & Z \ar[d] 
& \grl \ar@{}[l]|{\ni}  \ar[d]\\
(0,\grl,s) \ar@{}[r]|{\in} & \goM(d,v) \ar[r]^{\tilde p} & \goZ  
& (0,\grl) \ar@{}[l]|{\ni}  
}
$$
By Proposition \ref{prp:KNM2} the diagram is a pull back and by lemma
\ref{lem:lemma38} $p$ and $\tilde p$ are locally trivial over $R$, $\calR$.
We call $M_R = p^{-1}(R)$ and $\goM_{\calR}= \tilde p ^{-1} (\calR)$. 

Now consider the complex 
$\calF = R p_* (\mZ_{M_{R}}) = \mi_M^{-1} R\tilde p _* (\mZ_{\goM_{\calR}})  $ 
which is cohomologically a locally constant complex. 
We observe now that $\Pi_1(\calR)$ is trivial so $ R\tilde p _* (\mZ_{\goM_{\calR}})$
is isomorphic to cohomologically constant complex on $\calR$ so it is 
$\calF$ on $R$. In particular for any $x,y\in R$ we have a canonically isomorphism
$$
\psi^i_{x,y} : H^i(\calF_x) \lra H^i(\calF_y).
$$
Now observe that by Proposition \ref{prp:azioneweyl}
there is an action on $W$ on $R, M_R$ and that $p$ is equivariant with respect
to this action. So we can define a $W$ action on 
$H^i(M_{0,\grl}(d,v),\mZ)$ by
$$
 \grs (c) = \psi^i_{\grs\grl,\grl}\comp H^i(\Phi^{d,v}_{\grs,0,\grl}) (c) 
$$
for any $\grs \in W$.
To verify that this is an action we have only to verify that
$$
\psi^i_{\grs\grl^2,\grl^1}\comp H^i(\Phi^{d,v}_{\grs,0,\grl^2}) 
\psi^i_{\grl^1,\grl^2}(c)
=\psi^i_{\grs\grl^1,\grl^1}\comp H^i(\Phi^{d,v}_{\grs,0,\grl^1}) (c).
$$
Since $R$ is connected and $H^i(M,\mZ)$ is discrete this is clear.
S o we have proved the following corollary.
\begin{cor}
If $d=v$ and $(0,\grl) \in \calR$ then there is an action of $W$ on 
$H^i(M_{0,\grl}(d,v),\mZ)$.
\end{cor}

\begin{oss} If $m_+=(1,\dots,1)$ and $\grl =0$ it is easy to see 
that $d\grm_s$ is surjetive for all $s \in \grL_{m_+,0}(d,v)$. Then  by lemma 
\ref{lem:lemma38} there is a canonical isomorphism  
$H_*(M_{m_+,0}(d,v)) \isocan  H_*(M_{0,\grl}(d,v))$ if 
$(0,\grl) \in \calR$. So
by Nakajima's Theorem (Theorem 10.2 \cite{Na2}) it is natural to 
make the following conjecture:
\end{oss}
\begin{con} Let $top = \frac{1}{2} \dim H_*(M_{0,\grl}(d,v))$ then 
$$
H^{top}\bigl(M_{(0,\grl)}(d,v), \mC \bigr) \isocan \bigl(L_d \bigr)_0
$$
where $\bigl(L_d \bigr)_0$ is the $0$-weight space of the Kac-Moody algebra associated 
to the quiver of heighest weight $\sum_i d_i \bar {\omega}_i$.
\end{con}


\section{Reduction to the dominant case} \label{riduzioneM0}

As a consequence of Proposition \ref{prp:azioneweyl} we see
that if $(m,\grl) \in \calG_v$ then there exists $\grs \in W$ and 
$v' = \grs \cdot v$ such that $d-  v'$ is dominant  and 
$M_{\grs m , \grs \grl}(d,v') \isocan M_{m,\grl}(d,v)$.
We generalize now this result to arbitrary $\grl$.

On $Q$ we consider the following order:  $v' \leq v$
if and only if $v_i' \leq v_i$.

We consider now the following construction:
let $v' \leq v$ 
and fix an 
embedding $V'_i \incluso V_i$ and a complement $U_i$ of $V'$ in 
$V_i$, then we can define a map $\wt{\mj} \colon S(d,v') 
\lra S(d,v)$ through:
\begin{equation}\label{defembeddingj}
\wt{\mj}(B^{\prime},\grg^{\prime},\grd^{\prime})=
\left(
\begin{pmatrix}
B^{\prime} & 0 \\
0     &     0
\end{pmatrix}
,
\begin{pmatrix}
\grg^{\prime}  \\
0     
\end{pmatrix}
,
\begin{pmatrix}
\grd^{\prime} & 0 
\end{pmatrix}\right)
\end{equation}
where the matrices of the new triple represents the maps described 
through the decomposition $V_i=V'_i \oplus U_i$.

Suppose now that $(m_i,\grl_i)=0$ for all $i$ such that 
$v_i'\neq v_i$. Then  it is easy to see that this map restrict to a map  
$\mj_r : \grL_{m,\grl}(d,v') \lra \grL_{m,\grl}(d,v)$ and so enduces a map
$\mj^{v'}_v=\mj \colon M_{0,\grl}(d,v') \lra M_{0,\grl}(d,v)$.

\begin{lem} \label{lemmaj}
$\mj$ is a closed immersion
\end{lem}
\begin{proof}
We prove that the map $\mj^{\sharp} : \mC[\grL_{\grl}(d,v)]^{G(v)}
\lra \mC[\grL_{\grl}(d,v')]^{G(v' )}$ is surjective.
By proposition \ref{generatoriinvarianti} this follows 
by the following two identities:
$$
\Tr\left(\gra\left(\mj (s) \right) \right) = 
\Tr\left(\gra (s) \right) \; \mand \;
\grb\left(\mj (s) \right) = \grb (s) 
$$
for each $B$-path $\gra$ and for each admissible path $\grb$.
\end{proof}

\begin{lem} \label{lemmariduzioneM0}
If $2 v_i > d_i + \sum_{j \in I} a_{ij}v_j$ and $v' = v -\gra_i$
then $\mj$ is an isomorphism of algebraic varieties
\end{lem}
\begin{proof}
 It's enough to prove that $\mj$ is surjective. 
Let $s=(B,\grg,\grd) \in \grL_0(d,v)$ and consider the sequence (see \eqref{defTab}
for the notation) :
$$
\begin{CD}
T_i @>{b_i}>> V_i @>{a_i}>> T_i.
\end{CD}
$$
Since $b_ia_i =0 $ and $2 \dim V_i > \dim T_i$ we have that $b_i$ is not 
surjective or that $a_i$ is not injective. 

Suppose that $b_i$ is not surjective, then up to the action of $G_v$ we 
can assume that $\Im b_i \subset v'_i$. Then, for $t \in \mC^*$ consider 
$g_t=(g_{j,t}) \in G_v$ with 
$$
g_i =
\begin{pmatrix}
\Id_{v'_i} & 0 \\
0 & t^{-1}
\end{pmatrix}
\; \mand \;
g_j = \Id_{V_j} \mfor j \neq i.
$$
Then 
\begin{enumerate}
	\item $g_{i,t} B_h = B_h $ if $h_1 =i$ and $g_i \grg_i =\grg_i$, 
		since $\Im B_h , \Im \grg_i \subset \Im b_i \subset v'_i$,
	\item $\exists\, \lim_{t \to 0}  B_h g_{i,t}^{-1} = B_h $ if 
		$h_0 =i$ and $ \grd_i g_i^{-1} =\grd_i$
\end{enumerate}
So $\exists\, \lim_{t \to 0} g_t s = s^{\prime}$ and it is clear that  
$s^{\prime} \in \wt{\mj}(\grL_0(d,v'))$ and that $p_0(s)= p_0(s^{\prime})\in
\Im \mj$.

If $b_j$ is surjective and $a_i$ is not injective the argument is similar.
\end{proof}

\begin{prp}\label{prp:riduzioneM0}For all $\grl$ and for all $d\geq 0, v\geq 0$ 
there exists $v'$ and $\grs \in W$ such that $\bar d - \bar v ' $ is dominant and 
$$
M_{0 ,\grs \grl}(d,v') \isocan M_{0,\grl}(d,v).
$$
\end{prp}

\begin{proof}
We prove this proposition by induction on the order $\leq$ on $Q$.

\noindent{First step: $v=0$}. If $v=0$ we can take $v'=v$ and $\grs =1$.

\noindent{Inductive step.} If $d-v$ is not dominant then there exists 
$i$ such that $2v_i > d_i + \sum a_{ij}v_j$. 

If $ \grl_i \neq 0 $ we observe that $s_i v= v' < v$ (that is $v' \leq v$ 
and $v'\neq v$) and that $M_{s_i m,s_i \grl}(d,v')\isocan M_{m,\grl}(d,v)$
and so we can apply the  inductive hypothesis.

If $\grl_i = 0$ we apply the previous lemma and the inductive hypothesis.
\end{proof}


\section{On normality and connectdness of quiver variety in the finite type case}

In this section we resctrict our attention to the case
of quiver varieties of finite type and to the case
$m=(1,\dots,1)$ and $\grl =0$ and we fix $d,v$.
By remark \ref{oss:riduzioneanonzero}
we can assume without loss of generality that $v_i >0$ for all $i$.
We would like to prove the following conjecture:
\begin{con}
$M_{m,0}(d,v)$ is connected and $M_{0,0}(d,v)$ is normal. 
\end{con}

\begin{oss}
If $\tilde m = (m_1,\dots,m_n) \in\mN_+^I$ it is easy to see that 
$\grL_{\tilde m,0}= \grL_{m,0}$ ao in particular  $M_{m,0}(d,v)$ is smooth.
Instead $M_{0,0}$ is a cone so it is clearly connected.
\end{oss}

By proposition \ref{prp:riduzioneM0}it is enough to prove the theorem in the case
$d-v$ dominant. Unfortunately I'm not able to prove the conjecture only in the case
$d-v$ regular: $\bra d-v,\gra_i\ket >0$ for all $i$.
To prove the conjecture in this case we will use the following stratification introduced by
Lusztig in \cite{Lu:Q4}.

\begin{dfn}
For any $s\in S$ and $i\in I$
let 
$$
V_i^+ = V_i^+(s) = \sum_{\gra \text{ a }B-\text{path}\st \gra_1=i} \Im (\gra(s)\grg_{\gra_0})
$$
If $v'=(v_1',\dots,v_n') \in \mN^n$ we define
$$
\grL^{v'} = \{ s\in \grL_0(d,v) \st \dim V_i^+ (s) = v_i'\}.
$$
\end{dfn}

Observe that $\grL^v = \grL_{m_+,0}(d,v)$. To prove our result 
we will use the following lemma of Lusztig.

\begin{lem}[Lusztig: \cite{Lu:Q4} Proposition 4.5 and Proposition 5.3]
\label{lem:Lusztigdim}If $0\leq v'_i \leq v_i$ for each $i$ then 
$$
\dim \grL^{v'}(d,v) = 
\dim S - \sum_{i\in I} \dim gl(V_i) - \bra (v-v')\cech , d-v \ket - \tfrac{1}{2} 
	\bra (v-v')\cech , v-v' \ket
$$
\end{lem}

Our result follows trivially from the following lemma.

\begin{lem}
1) If $d-v$ is dominant then $\grL_{0,0}(d,v)$ is a complete intersection.

\noindent 2) If $d-v$ is regular then $\grL_{0,0}$ is normal and
  		irreducible and $\grL_{m_+,0}(d,v)$ is connected.
\end{lem}
\begin{proof}
Observe that $\grL_{0,0}(d,v) =\mu ^{-1}(0)$ 
so each irreducible component of $\grL_{0,0}(d,v)$ must have 
dimension at least $\dim S - \sum_i \dim gl(V_i) = \delta_V$.

Suppose now that $d-v$ is dominant. By Nakajima's theorem (\cite{Na2} Theorem
10.2)  $M_{m_+,0}$ is not empty.
Observe also that
by Proposition \ref{prp:1.14} $\grL_{m_+,0}(d,v)$ is a smooth
subset of $\grL_{0,0}(d,v)$
of dimension $\delta_V$.

It is well known that $\grL_{m_+,0} (d,v) = \grL^v$. Hence 
$$
\grL_{0,0}(d,v) - \grL_{m_+,0} (d,v)  = \bigcup_{v'\leq v \mand v' \neq v} \grL^{v'}.
$$
By the lemma above  we have that if $v' \leq v$ and $v' \neq v$
$\dim  \grL^{v'} < \delta_V$. So $\grL_{m_+,0}(d,v)$ must be dense in $\grL_{0,0}(d,v)$ 
and  $\grL_{0,0}(d,v)$ is a complete intersection. Moreover if $d-v$ is regular
we have that $\dim  \grL^{v'} < \delta_V - 1$ so the singular locus has codimension 
at least two and normality and irreducibility follows. 
Finally by our discussion it is clear that if $\grL_{m_+,0}(d,v)$ is disconnected
then $\grL_{0,0}(d,v)$ is not irreducible.
\end{proof}

\begin{oss}
In the lemma we can substitute $\grL_{m_+,0}(d,v)$ 
with any other subset $Reg$ of regular points in $\grL(d,v)$.
In this way is indeed possible to improve a little bit the theorem but
Crawley-Boevey explained me that this strategy cannot work in general becouse
there are cases where $d-v$ is dominant and $\grL_{0,0}(d,v)$ is not normal.
It should be also pointed out that Crawley-Boevey proved the connectdness
in complete generality (\cite{CrBo}). He said me that is also 
able to prove normality
for a much bigger class of quiver varieties.
\end{oss}


\bibliographystyle{amsplain} 
\bibliography{tato,quiver,rapprese}




\end{document}